\documentclass[reqno,a4paper,11pt]{article}

\usepackage{amsmath,amssymb,amsthm}

\usepackage{mathrsfs}

\usepackage{mathtools}
\usepackage{commath}

\usepackage{thmtools}
\usepackage{thm-restate}

\usepackage{cases}
\usepackage{enumitem}
\setlist[enumerate,1]{label=(\roman*)}

\usepackage{xcolor}
\definecolor{softred}{RGB}{220,25,25}
\usepackage[pdftex, pdfborderstyle={/S/U/W 0}]{hyperref}
\hypersetup{
    colorlinks=true,
    linkcolor=softred,
    citecolor=blue,
	hypertexnames=false,
	urlcolor=red,
}

\usepackage{cleveref}

\usepackage{etoolbox}

\usepackage{comment}

\usepackage[numbers]{natbib}

\numberwithin{equation}{section}

\usepackage[framemethod=default]{mdframed}
\usepackage{marginnote}
\usepackage{environ} 
\usepackage{geometry}

\makeatletter

\newif\ifthmboxes
\thmboxesfalse

\newif\ifproofbar
\proofbarfalse

\newif\iframbling
\ramblingfalse

\colorlet{clr@theorem}{red}
\colorlet{clr@proposition}{red}
\colorlet{clr@lemma}{red}
\colorlet{clr@corollary}{red}
\colorlet{clr@claim}{red}
\colorlet{clr@fact}{red}
\colorlet{clr@observation}{red} 

\colorlet{clr@question}{magenta}
\colorlet{clr@conjecture}{magenta}
\colorlet{clr@openproblem}{magenta}

\colorlet{clr@definition}{green}
\colorlet{clr@construction}{green}

\colorlet{clr@example}{blue}

\definecolor{clr@remark}{HTML}{9D4EDD} 

\colorlet{clr@proof}{black!30}   

\definecolor{clr@comment}{HTML}{2A9D8F}   

\newcommand{\marginSideNormal}{25mm}
\newcommand{\marginSideRambling}{45mm}
\newcommand{\marginNoteWidth}{35mm}
\newcommand{\marginNoteSep}{5mm}
\newcommand{\marginTopNormal}{25mm}
\newcommand{\marginBottomNormal}{25mm}
\iframbling
\else
\fi

\newcommand{\thmboxsetup}[1]{%
  linecolor=#1,%
  linewidth=0.6pt,%
  backgroundcolor=#1!8,%
  roundcorner=2pt,%
  innertopmargin=6pt,%
  innerbottommargin=6pt,%
  innerleftmargin=8pt,%
  innerrightmargin=8pt,%
  skipabove=8pt,%
  skipbelow=8pt,%
}

\newcommand{\newthmboxstyle}[3]{%
  \ifthmboxes
    \begingroup
    \edef\thmbox@opts{\thmboxsetup{#2}}%
    \expandafter\endgroup
    \expandafter\newthmboxstyle@i\expandafter{\thmbox@opts}{#1}{#3}%
  \else
    \declaretheoremstyle[
      headfont=\normalfont\bfseries,
      bodyfont=#3,
    ]{#1}%
  \fi
}
\newcommand{\newthmboxstyle@i}[3]{%
  \declaretheoremstyle[
    headfont=\normalfont\bfseries,
    bodyfont=#3,
    mdframed={#1},
  ]{#2}%
}

\newcommand{\newthmhiddenstyle}[5]{%
  \ifthmboxes
    \begingroup
    \edef\thmbox@opts{\thmboxsetup{#2}}%
    \expandafter\endgroup
    \expandafter\newthmboxstyle@i\expandafter{\thmbox@opts}{#1}{#3}%
    \declaretheorem[style=#1,numbered=no,name=#5]{#4}%
  \else
    \NewEnviron{#4}{{#3\ignorespaces\BODY}}%
  \fi
}

\newthmboxstyle{box-theorem}{clr@theorem}{\itshape}
\newthmboxstyle{box-proposition}{clr@proposition}{\itshape}
\newthmboxstyle{box-lemma}{clr@lemma}{\itshape}
\newthmboxstyle{box-corollary}{clr@corollary}{\itshape}
\newthmboxstyle{box-claim}{clr@claim}{\itshape}
\newthmboxstyle{box-fact}{clr@fact}{\itshape}
\newthmboxstyle{box-observation}{clr@observation}{\itshape}
\newthmboxstyle{box-question}{clr@question}{\itshape}
\newthmboxstyle{box-conjecture}{clr@conjecture}{\itshape}
\newthmboxstyle{box-openproblem}{clr@openproblem}{\itshape}
\newthmboxstyle{box-definition}{clr@definition}{\normalfont}
\newthmboxstyle{box-construction}{clr@construction}{\normalfont}
\newthmboxstyle{box-example}{clr@example}{\normalfont}
\newthmboxstyle{box-remark}{clr@remark}{\normalfont}

\makeatother

\ifproofbar
  \surroundwithmdframed[
    topline=false,
    bottomline=false,
    rightline=false,
    leftline=true,
    linewidth=1.2pt,
    linecolor=clr@proof,
    innertopmargin=2pt,
    innerbottommargin=2pt,
    innerleftmargin=10pt,
    innerrightmargin=0pt,
    skipabove=6pt,
    skipbelow=6pt,
  ]{proof}
\fi

\iframbling
  \NewEnviron{rambling}{%
    \begin{mdframed}[
        linecolor=clr@comment,
        backgroundcolor=clr@comment!8,
        topline=false, bottomline=false, rightline=false, leftline=true,
        linewidth=2pt,
        roundcorner=0pt,
        innertopmargin=5pt, innerbottommargin=5pt,
        innerleftmargin=8pt, innerrightmargin=8pt,
        skipabove=6pt, skipbelow=6pt,
      ]%
      \sffamily\small
      \textcolor{clr@comment}{\bfseries Commentary.}\hspace{0.5em}%
      \ignorespaces\BODY
    \end{mdframed}%
  }
\else
  \NewEnviron{rambling}{}
\fi

\declaretheorem[style=box-theorem,numberwithin=section,name=Theorem]{theorem}
\declaretheorem[style=box-proposition,sibling=theorem,name=Proposition]{proposition}
\declaretheorem[style=box-lemma,sibling=theorem,name=Lemma]{lemma}

\declaretheorem[style=box-conjecture,sibling=theorem,name=Conjecture]{conjecture}

\declaretheorem[style=box-question,sibling=theorem,name=Question]{question}

\declaretheorem[style=box-theorem,numbered=no,name=Theorem]{theorem-n}
\declaretheorem[style=box-proposition,numbered=no,name=Proposition]{proposition-n}
\declaretheorem[style=box-lemma,numbered=no,name=Lemma]{lemma-n}
\declaretheorem[style=box-corollary,numbered=no,name=Corollary]{corollary-n}
\declaretheorem[style=box-conjecture,numbered=no,name=Conjecture]{conjecture-n}
\declaretheorem[style=box-claim,numbered=no,name=Claim]{claim-n}
\declaretheorem[style=box-fact,numbered=no,name=Fact]{fact-n}
\declaretheorem[style=box-openproblem,numbered=no,name=Open Problem]{openproblem-n}
\declaretheorem[style=box-question,numbered=no,name=Question]{question-n}
\declaretheorem[style=box-observation,numbered=no,name=Observation]{observation-n}

\declaretheorem[style=box-remark,sibling=theorem,name=Remark]{remark}
\declaretheorem[style=box-definition,sibling=theorem,name=Definition]{definition}

\newthmhiddenstyle{hidden-remark}{clr@remark}{\normalfont}{remark-n}{Remark}
\newthmhiddenstyle{hidden-definition}{clr@definition}{\normalfont}{definition-n}{Definition}
\newthmhiddenstyle{hidden-construction}{clr@construction}{\normalfont}{construction-n}{Construction}
\newthmhiddenstyle{hidden-example}{clr@example}{\normalfont}{example-n}{Example}

\newcommand{\defined}{\mathrel{\coloneqq}}

\DeclarePairedDelimiter{\p}{\lparen}{\rparen}

\let\b\relax
\DeclarePairedDelimiter{\b}{\lbrack}{\rbrack}

\newcommand{\st}{\mathbin{\colon}}

\undef{\set}
\DeclarePairedDelimiter{\set}{\lbrace}{\rbrace}

\undef{\emptyset}
\newcommand{\emptyset}{\varnothing}

\newcommand{\union}{\mathbin{\cup}}
\newcommand{\inter}{\mathbin{\cap}}

\DeclareMathOperator{\ind}{\mathbf{1}}

\newcommand{\from}{\colon}

\let\d\relax
\newcommand{\d}{\, \mathrm{d}}

\DeclareMathOperator{\loglog}{log\,log}

\undef{\mod}
\newcommand{\mod}[1]{\ (\mathrm{mod}\ #1)}

\undef{\abs}
\DeclarePairedDelimiterX{\abs}[1]
  {\lvert}{\rvert}{\ifblank{#1}{\,\cdot\,}{#1}}

\undef{\norm}
\DeclarePairedDelimiterX{\norm}[1]
  {\lVert}{\rVert}{\ifblank{#1}{\,\cdot\,}{#1}}

\DeclarePairedDelimiterX{\inner}[2]
  {\langle}{\rangle}{\ifblank{#1}{\,\cdot\,}{#1},\ifblank{#2}{\,\cdot\,}{#2}}

\DeclarePairedDelimiterX{\absinner}[2]
  {|\langle}{\rangle|}{\ifblank{#1}{\,\cdot\,}{#1},\ifblank{#2}{\,\cdot\,}{#2}}

\DeclareMathDelimiter{\given}
  {\mathbin}{symbols}{"6A}{largesymbols}{"0C}

\DeclareMathOperator{\Prob}{\mathbb{P}}
\DeclarePairedDelimiterXPP{\prob}[1]
  {\Prob}{\lparen}{\rparen}{}
  {\renewcommand{\given}{\nonscript\;\delimsize\vert\nonscript\;\mathopen{}}#1}

\DeclareMathOperator{\Expec}{\mathbb{E}}
\DeclarePairedDelimiterXPP{\expec}[1]
  {\Expec}{\lparen}{\rparen}{}
  {\renewcommand{\given}{\nonscript\;\delimsize\vert\nonscript\;\mathopen{}}#1}

\DeclareMathOperator{\Var}{Var}
\DeclarePairedDelimiterXPP{\var}[1]
  {\Var}{\lparen}{\rparen}{}
  {\renewcommand{\given}{\nonscript\;\delimsize\vert\nonscript\;\mathopen{}}#1}

\DeclareMathOperator{\Cov}{Cov}
\DeclarePairedDelimiterXPP{\cov}[2]
  {\Cov}{\lparen}{\rparen}{}{#1,#2}

\newcommand{\eps}{\varepsilon}
\newcommand{\sseq}{\subseteq}

\let\l\relax
\newcommand{\l}{\ell}

\let\SS\relax

\newcommand{\CC}{\mathbb{C}}

\newcommand{\EE}{\mathbb{E}}

\newcommand{\RR}{\mathbb{R}}
\newcommand{\SS}{\mathbb{S}}

\newcommand{\zz}{\zeta}

\usepackage{graphicx,tikz}
\usetikzlibrary{arrows.meta,patterns,shapes.misc}

\usepackage{scrextend}

\usepackage[T1]{fontenc}
\usepackage{titlesec}

\titleformat{\section}{\centering\bfseries\scshape\Large}{\thesection}{1em}{}
\titleformat{\subsection}{\bfseries\scshape\large}{\thesubsection}{1em}{}

\newcommand{\blank}{\;\cdot\;}

\newcommand{\signs}{\set{-1,+1}}
\newcommand{\cube}{\set{0,1}}
\newcommand{\lth}{$\l$\textsuperscript{th}}
\newcommand{\kth}{$k$\textsuperscript{th}}
\newcommand{\upstairs}[2]{U_{#1}^{#2}}
\newcommand{\downstairs}[2]{L_{#1}^{#2}}

\newcommand{\dTV}{d_{\text{TV}}}

\begin{document}

\title{\textsc{\bfseries Proving the Koml\'os conjecture slowly}}

\renewcommand{\thefootnote}{\fnsymbol{footnote}}

\author{\textsc{Lawrence Hollom}\footnotemark[1]}

\footnotetext[1]{\href{mailto:lawrence.hollom@epfl.ch}{lawrence.hollom@epfl.ch}. Institute of Mathematics, EPFL, Switzerland.}

\renewcommand{\thefootnote}{\arabic{footnote}}

\date{}

\maketitle

\begin{abstract}
  The Koml\'os conjecture was a long-standing open problem in discrepancy theory until it was recently solved in a flurry of AI developments.
  We give a motivated exposition of an elementary proof of the conjecture due to Karingula and Lovett, as well as a counterexample to the strong Koml\'os conjecture due to Kintali, aiming to introduce the ideas in a gentle and explained manner in order to ease the understanding of the proofs.
  We conclude with several problems which remain open.
\end{abstract}

\begingroup
\hypersetup{linkcolor=black}
\tableofcontents
\endgroup
\newpage


\section{Introduction}
\label{sec:intro}

The Koml\'os conjecture asks whether a collection of vectors can be split into two balanced groups.
More precisely, it asserts that such a split is possible: given $v_1,\dotsc,v_n \in \RR^d$, assign each vector $v_i$ a sign $\eps_i \in\signs$ so that the vector $\sum_i \eps_i v_i$ has no large coordinate in absolute value.
The conjecture is formally stated as follows.

\begin{conjecture}[The Koml\'os Conjecture]
	\label{conj:komlos}
	There is a constant $C > 0$ such that, for all sequences $v_1,\dotsc,v_n\in\RR^d$ with $\norm{v_i}_2 \leq 1$, there are $\eps_1,\dotsc,\eps_n \in \signs$ such that
	\begin{equation*}
		\norm[\bigg]{\sum_{i=1}^n \eps_i v_i}_\infty \leq C.
	\end{equation*}
\end{conjecture}

A proof that \Cref{conj:komlos} holds was recently found by an AI system, prompted by Guo, Fang, and Lu \cite{GFL26a} and building on ideas of Smirnov and Vershynin \cite{SV26}, among others. (See \Cref{sec:history} for a thorough history of the problem.)
The value of $C$ achieved in \cite{GFL26a} was $C = 3\sqrt{2\pi}$.
We give a different proof, primarily derived from that of Karingula and Lovett \cite{KL26}, which we believe to be easier to follow.
In particular, we prove the following.

\begin{theorem}
	\label{thm:main}
	\Cref{conj:komlos} holds with $C = 18\sqrt{3}$.
\end{theorem}

The second result we present is a proof due to Kintali \cite{Kin26} that the ``strong Koml\'os conjecture'' is false.
More precisely, the result is as follows.

\begin{theorem}
	\label{thm:false}
	There is a constant $c > 0$ such that, for all positive integers $d$ there is an integer $n = \exp(O(d^2))$ and vectors $v_1,\dotsc,v_n\in\RR^d$ with $\norm{v_i}_2 = 1$ for all $i$ such that, for all signings $\eps_1,\dotsc,\eps_n\in\signs$, there is an integer $1\leq k\leq n$ such that
	\begin{equation*}
		\norm[\bigg]{\sum_{i=1}^k \eps_i v_i}_\infty \geq c\sqrt{\log d}.
	\end{equation*}
\end{theorem}

Although we will not use this fact, note that one may assume that $n=d$ in \Cref{conj:komlos}: if we are given $n\leq d$ vectors, then we can add copies of the zero vector until $n=d$, which does not affect the signed sum of the vectors.
If $n\geq d$, then we can append zero entries to all the vectors until $d = n$.
In this way, we can consider the Koml\'os conjecture as concerning a square $n\times n$ matrix $A$ with normalised columns; the conjecture asks whether there is a vector $\eps\in\signs^n$ with $\norm{A\eps}_\infty \leq C$.

Thinking geometrically, $\norm{x}_\infty \leq C$ is equivalent to $x$ lying inside the axis-aligned hypercube of side length $2C$ centred at the origin.
This leads to the following equivalent phrasing of \Cref{conj:komlos}:
is there a linear transformation $A$ from $\RR^n$ to $\RR^n$ such that $\norm{Ae_i}_2 \leq 1$ for every basis vector $e_i$ which maps every vertex of the hypercube $\signs^n$ to outside of the larger hypercube $[-C,C]^n$?
The conjecture is that there should be no such transformation.
One particular case of this is that of a rotation: can we rotate a small hypercube so all its vertices lie outside of a big hypercube?
While it is possible to rotate many points of $\signs^n$ to be outside of $[-C,C]^n$, there are an exponential number of points which must be dealt with simultaneously.
In this light, the conjecture seems plausible.
However, high-dimensional geometry can be notoriously unintuitive, and so this explanation alone is not necessarily a good reason to believe the conjecture.


\subsection{History and related problems}
\label{sec:history}

We now provide a thorough account of the history of both this problem and surrounding conjectures. 
The narrative is broken into two parts: one covering the progress made on the conjecture before a proof was found, and another on the flurry of activity in the weeks after the release of the first proof.


\subsubsection{Before proof announcement}
\label{subsec:bc}

After being posed by Koml\'os in the early 1980s, \Cref{conj:komlos} was circulated by word-of-mouth in the mathematical community (see B\'ar\'any and Grinberg \cite[end of Section 2]{BG81} for the first written reference to the problem we could find).
The first significant partial result appeared in Spencer's seminal 1985 paper ``six standard deviations suffice'' \cite{Spe85}, which was also the origin of the name ``the Koml\'os conjecture''.
Spencer's results were independently discovered (up to the constants) on the other side of the iron curtain by Gluskin \cite{Glu89}.
Alongside bringing the conjecture to a wider audience, Spencer proved the first significant bound on the problem, proving that there are signs $\eps_1,\dotsc,\eps_n\in\signs$ for which $\norm{\sum_i \eps_i v_i}_\infty \leq O(\log d)$. 
The existence of some bound which depends only on the dimension had been known since Spencer's 1977 work on vector balancing games \cite{Spe77}.

The method behind Spencer's proof was the partial colouring technique pioneered by Beck \cite{Bec81}.
To say only a few words about this technique, to find a good choice of signs $\eps\in\signs^n$, one tracks a vector $x\in[-1,+1]^n$, varying the entries, locking in place any which hit the endpoints $\pm 1$.
By carefully updating this vector in steps, making sure that a constant proportion of the remaining coordinates of $x$ hit $\pm 1$ every step, one can reach a vertex in $\signs^n$ while accumulating little discrepancy.

Also preceding Spencer's work and using a similar technique was the paper of Beck and Fiala \cite{BF81} which gave rise to the eponymous Beck--Fiala conjecture, which is a special case of the Koml\'os conjecture.
This conjecture asserts that, given a family $S_1,\dotsc,S_n$ of subsets of a universe of size $d$ (which we label $[d] = \set{1,\dotsc,d}$) such that each element of $[d]$ is contained in at most $t$ sets $S_i$, there is an assignment $\chi \from [d] \to\signs$ such that the sum of these signs over each set, $\chi(S_i)$, satisfies $\abs{\chi(S_i)}\leq O(\sqrt{t})$.
To see that this is a special case of the Koml\'os conjecture, take the vectors $v_i$ to be indicator functions of which sets $S_j$ the element $i$ is contained in, divided by $\sqrt{t}$ (so that they are sub-units), i.e. $(v_i)_j = t^{-1/2} \ind_{i\in S_j}$.

\medskip

This area is rich with problems, and many variations of these questions have been considered: Spencer \cite{Spe86} also worked on the question in which the incoming vectors $v_i$ have $\norm{v_i}_\infty \leq 1$ (as opposed to the bound $\norm{v_i}_2 \leq 1$ which is our primary focus) and proved that, if $n=d$, then signs may be chosen for the vectors so that the infinity norm of the signed sum is at most $O(\sqrt{n})$.
This can be seen to be optimal by choosing all entries of all vectors $v_i$ independently and uniformly at random from $\signs$.
Moreover, he proved that one can enforce the same upper bound on all partial sums, but only if one has prior knowledge of the entire sequence of vectors.

The Koml\'os conjecture itself was then picked up by Hajela \cite{Haj88}, who believed it should be false.
He proved a lower bound which looked to be very close to achieving this (though we know now that the Koml\'os conjecture is true) and conjectured that there should be a sequence of vectors with all signed sums having infinity norm at least $\Omega(\sqrt{\log d})$. 
The intuition for this is that, given a uniform random orthogonal matrix $U\in \RR^{d\times d}$ and a fixed signing $\eps\in \signs^d$, the vector $U \eps$ is uniformly randomly distributed on the sphere of radius $\sqrt{d}$, and so the coordinates behave like standard Gaussians: the largest coordinate of such a vector is thus likely to be of order $\sqrt{\log d}$, leading to the belief that this should be the correct order of magnitude. 
This highlights a key tension in the problem: it is easy enough to force some signed sums---even to force many signed sums---to have large discrepancy, but there is an exponentially large number of signed sums to consider.

Work continued throughout the 90s on the Koml\'os and Beck--Fiala conjectures, with a number of related results emerging: 
Banaszczyk \cite{Ban90} proved an optimal Beck--Fiala-type result for the Euclidean norm;
Akemann and Anderson \cite{AA95} proved that the known bounds on the Beck--Fiala problem were tight for a certain continuous relaxation of the problem;
Bednarchak and Helm \cite{BH97} proved a small improvement on the best-known bound for the original Beck--Fiala problem;
Srinivasan discovered a slight improvement of Spencer's bound on the Koml\'os problem \cite{Sri97};
and Matou\v{s}ek \cite{Mat98} considered the (suitably normalised) problem in the $L_p$ norm, and proved the optimal result for all finite $p$.

Returning our narrative to the Koml\'os problem, in 1997 Giannopoulos \cite{Gia97} proved that, if $K\sseq \RR^d$ is a symmetric convex body of Gaussian measure $1/2$ and $v_1,\dotsc,v_d \in \RR^d$ are sub-unit vectors, then there is a signed sum $\sum_i \eps_i v_i$ which lies in $(c \log d) K$.
Banaszczyk \cite{Ban98} built on these results to achieve a significant strengthening of Spencer's bound on the Koml\'os problem: he proved that for all sequences $v_1,\dotsc,v_n \in \RR^d$ of sub-unit vectors, there is a signed sum satisfying $\norm[\big]{\sum_i \eps_i v_i}_\infty \leq O(\sqrt{\log n})$.
The key to Banaszczyk's bound was a result showing that, if a convex set $K \sseq \RR^d$ has Gaussian measure at least $1/2$ and $v$ is a sub-unit vector, then there is a convex set contained in $(K - cv) \union (K + cv)$ with Gaussian measure at least $1/2$ for some constant $c > 0$ ($c=1/5$ suffices).
Such a set must necessarily contain the origin if it is closed, and so, if $K$ is a hypercube of Gaussian measure at least $1/2$, there is always a point $x\in c^{-1}K$ and a choice $\eps\in\signs^n$ such that $x + \sum_i \eps_i v_i = 0$.
The bound of the form $\sqrt{\log d}$ comes from the required size of a hypercube to have Gaussian measure $1/2$.

Banaszczyk \cite{Ban12} would later revisit a related problem: can one balance not only the entire signed sum, but all prefixes? 
That is, given vectors $v_1,\dotsc,v_n$ with $\norm{v_i}_2 \leq 1$, choose signs to minimise $\max_k \norm{\sum_{i=1}^k \eps_i v_i}_p$.
For $p=2$, he proved an upper bound of $O(\sqrt{d} + \sqrt{\log n})$, and for $p=\infty$ (the setting of the Koml\'os conjecture), he proved $O(\sqrt{\log d + \log n})$; we will revisit this bound in \Cref{sec:strong-counter}.
Finally, he proved the optimal bound of $O(\sqrt{d})$ for prefix sums in the Euclidean norm if one can not only choose the signs but also a permutation of the order of the vectors.
We will return to this latter problem in \Cref{sec:conc}.

\medskip

We now pick up the thread of a related direction of research: so far, the results we have described for the Koml\'os problem are non-constructive, and do not lead to an efficient algorithm to find a signed sum of small discrepancy.
(Of course, an exponential time algorithm exists: just try every signing! By an ``efficient'' algorithm, we mean polynomial-time.)
A natural question is thus to ask not only for a bound on the signed sum, but a polynomial-time algorithm to construct a sequence of signs achieving this bound.
This question was considered by Bansal \cite{Ban10} in 2010, and he produced an efficient randomised polynomial-time algorithm to recover Spencer's $O(\log d)$ bound.
This introduced the technique to make Beck's technique \cite{Bec81} of partial colourings algorithmic by making small updates to variables $x_i \in [-1,+1]$, keeping coordinates fixed as soon as they hit $\pm 1$, and controlling the evolution of the process via a semidefinite program (SDP).
This technique was later refined and simplified by Lovett and Meka \cite{LM15}, in particular, the proof of Lovett and Meka was the first purely algorithmic proof, as Bansal's proof that the algorithm worked relied on the earlier non-algorithmic proof.

Also around this time in the mid 2010s, Nikolov \cite{Nik13} proved that the Koml\'os conjecture holds if one uses unit vectors in $\RR^n$ in place of the $\pm 1$ signs, and Bukh \cite{Buk16} improved the bound on the Beck--Fiala problem to $2t - \log^* t$.

The next improvement to the algorithmic bound came in 2016, when Bansal, Dadush, and Garg \cite{BDG16} found an efficient randomised algorithm matching Banaszczyk's $O(\sqrt{\log d})$ bound, refining the SDP approach to more carefully track how the fractional signing evolves in order to avoid particular adversarial conditions.
This approach was then adapted to other vector-balancing problems by Bansal and Garg \cite{BG17}.

A related line of work asked for these algorithmic approaches, previously randomised, to be made deterministic.
Dadush, Garg, Lovett, and Nikolov \cite{DGLN19} made a step in this direction, and then later Levy, Ramadas, and Rothvoss \cite{LRR17} found a fully deterministic algorithm recovering Banaszczyk's bound for the Koml\'os conjecture.
Finally, Bansal, Dadush, Garg, and Lovett \cite{BDGL17} extended the approach of \cite{DGLN19} to the full theory of Gaussian-convex bodies.

We now reach the 2020s.
Alweiss, Liu, and Sawhney \cite{ALS21} gave a linear-time randomised algorithm (faster than previous known algorithms) for achieving $O(\log (dn))$ discrepancy.
With work on the upper bounds continuing, Kunisky \cite{Kun21} gave the best-known lower bound on the constant in the Koml\'os conjecture: $1 + \sqrt{2}$, which still stands as the largest such bound.

Interested in the robustness of the conjecture, Bansal, Jiang, Meka, Singla, and Sinha \cite{BJMSS22a,BJMSS22b} considered a ``smoothed'' version of the conjecture, wherein a small, independent perturbation with suitably far from degenerate covariance is added to each vector.
They proved two main results in this direction.
First, in \cite{BJMSS22a} they proved an upper bound of $O(\sqrt{\log d + \loglog n})$ for the Koml\'os prefix problem in the above setting of added noise.
They also used the name ``strong Koml\'os conjecture'' to refer to this prefix problem, wherein one must bound not only the signed sum of the vectors, but $\norm{\sum_{i=1}^k \eps_i v_i}_\infty$ for all $k\leq n$.
They also proved \cite{BJMSS22b} the Koml\'os conjecture itself holds with high probability in the context of perturbing by Gaussian noise, provided that $n = \omega(d \log d)$.

\medskip

Finally, after 27 years, there was an improvement to the bound on the original Koml\'os problem, beating Banaszczyk's $O(\sqrt{\log d})$: Bansal and Jiang \cite{BJ25} proved $\widetilde{O}(\log^{1/4} d)$, where the tilde indicates hidden polyloglog factors.
Their techniques also yielded strong partial results on the Beck--Fiala conjecture, proving it true for a large range of parameters.
An exposition of this result with a simpler proof was recently given by the same authors \cite{BJ26}.

Kulkarni, Reis, and Rothvoss \cite{KRR24} proved sharp bounds for an online version of the Koml\'os conjecture, where the vectors arrive one at a time with no foresight of future vectors, improving on the results of \cite{ALS21}.

Now reaching 2026, we approach the conclusion of this first section of history.
Ercan \cite{Erc26} removed the $\log\log$ factors from the proof of Bansal and Jiang, proving an upper bound of $O(\log^{1/4} d)$ for the Koml\'os conjecture.
In a similar vein to Nikolov's work \cite{Nik13} relaxing the Koml\'os conjecture to vector colourings, Guillen and Kobzar \cite{GK26} proved that one can relax the condition of $\eps_i\in\signs$ to $\eps_i\in\CC$ being a unit-modulus complex number and the conjecture does hold.
Altschuler and Tikhomirov \cite{AT26} proved optimal bounds for the online Beck--Fiala problem.

Finally, we mention a result which may at first seem only tangentially related: Smirnov and Vershynin \cite{SV26} considered random walks of the form $S_k = \sum_{i=1}^k \xi_i v_i$, where $\xi_i$ are independent Rademacher random variables (i.e.\ distributed uniformly on $\set{-1,+1}$) and asked whether all $S_k$ can be kept within some symmetric convex set $K\sseq \RR^d$ by discarding a small proportion of the steps. 
In particular, they proved that this was possible, relating the number of vectors which were discarded to the Fisher information.
As discussed in \cite{GFL26a}, the ideas of Smirnov and Vershynin were highly influential in the first proof of the Koml\'os conjecture.

This brings us near to the present day (at time of writing), and the solution of the conjecture.


\subsubsection{Proof announcement onwards}
\label{subsec:ad}

We now discuss the results which have emerged since a proof of \Cref{thm:main} was discovered.
The release of the papers discussed in this section has been so rapid that there is little use talking about the month of their release, let alone the year: from start to end this sub-section covers only thirty days.
For this reason we do not present the papers in global chronological order, but rather grouped by how they attack the problem: there are three main proofs to discuss here: the method of ``low-energy witness functions'' which led to the first proof, the elementary ``collapse'' proof (which we present an exposition of), and the algorithmic proof.

Kicking off the flurry of activity, Guo, Fang, and Lu \cite{GFL26a} released a proof (discovered by the ``Odin Automatic AI Research Agent'') of \Cref{thm:main}.
The proof was quite technical, making use of significant machinery.
The constant bound given by this approach was improved from approximately $7.5$ to approximately $6.9$ by a further AI-assisted proof by Akbas and Sra \cite{AS26b}.
Akbas and Sra also found a more general version of this result \cite{AS26a}, which is a matrix small-ball inequality, which implies several prominent results from the area such as the Koml\'os conjecture, and matrix Spencer.
It turned out that much of the heavy machinery for this proof is not in fact required if one is not seeking to optimise the constant, and an exposition of a simplified version of this proof due to Bandeira can be found at the Randomstrasse 101 blog \cite{Ban26a}.

Moving further in the direction of removing heavy machinery, a proof of the Koml\'os conjecture was found by Karingula and Lovett \cite{KL26}, and this is the version of the proof that our exposition is based on.

Guo, Fang, and Lu \cite{GFL26b} further used their AI system to find an efficient algorithm giving constant discrepancy; a version with faster running time was subsequently reported by Li \cite{Li26}.
An exposition of a simplified version of this algorithmic proof, highlighting how the approach differed from previous attempts, has been produced by Bansal \cite{Ban26b}.

The final result we mention here is the proof of \Cref{thm:false}, due to Kintali \cite{Kin26}, showing that one cannot hope for a constant bound if one wishes to control the infinity norm of all partial sums.
We present a version of this proof in \Cref{sec:strong-counter}.


\subsection{Proof outline}

In the interest of making the proof of \Cref{thm:main} as easy to understand as possible, we will go through it three times at varying levels of detail.
First, in the present section, we outline the proof at a high level, which should give a sense of how the overarching strategy works. 
Then, in \Cref{subsec:collapse-motivation}, we provide a motivation for why the ideas and tools we use are reasonable, and how one might arrive at them (which the author has attempted to make sufficiently thorough that one could fill in all the details without looking at the proofs). 
Finally, in \Cref{sec:back}, we give the full formal details.

The proof we present here is completely elementary and does not require any of the spectral machinery required by other proofs of the conjecture (see e.g.\ Bandeira's post on the Randomstrasse 101 blog \cite{Ban26a} and Bansal's expository paper \cite{Ban26b} for excellent presentations of two different proofs using spectral machinery).
The downside, however, is that this proof is perhaps more difficult to motivate and describe, and takes place in the slightly strange-looking space $\RR^d \times \cube^n$.

\medskip

The goal of this proof is to exhibit a probability distribution supported on $K_0 \defined (-C,C)^d$ with mean $\sum_i \eps_i v_i$ for some $\eps_i = \pm 1$, which demonstrates that this signed sum is in $K_0$, as required.
For a probability density $P \from K \to \RR$, we write $\mu(P) \defined \int_K x \, P(x) \d x$ for its centre of mass, or mean.
We will adjust this mean by one vector at a time.
To move $\mu(P)$ by $\eps_i v_i$, we constrain ourselves to only move mass from points $x$ to points $x-av_i$ or $x+av_i$ for some constant $a > 0$ (in fact, $a = 5$ will suffice). 
This way, we only need to be able to move a sufficient constant fraction of the mass in one direction in order to move the mean by $\pm v_i$.
As such an operation can distort the measure $P$ significantly, and we wish to perform $n$ such operations in a row, we carefully prepare the ground, first moving to a larger space which allows us to guarantee that we can still move mass in the desired way.

We start with a density $P_0$ on $K_0$ with $\mu(P_0) = 0$ and ``lift'' it to successively larger spaces $K_k \sseq \RR^d \times \cube^k$.
Assume that the two densities $P_0 - av_1$ and $P_0 + av_1$ are similar, in that the total variation distance between them is small.
This means that there is a lot of mass at points $x$ such that both $x-av_1$ and $x+av_1$ are in $K_0$; we reserve some of this mass and move it to a separate set: this separate set is contained in $\RR^d \times\set{1}$, whereas the rest of the mass stays in $\RR^d\times\set{0}$, giving us a measure $P_1$ on $\RR^d\times\cube$.
We will also be able to enforce that $\mu(P_1) = \mu(P_0) = 0$.
If we perform this lift operation in the right way, then $P_1$ has the same total variation measure property as $P_0$, and so we can iterate this lifting procedure: after $k$ steps, we have a density $P_k$ on $K_k\sseq \RR^d \times \cube^k$.
We emphasise here that this is not a collection of $2^k$ separate densities, but one single density (the total mass over all the $2^k$ parts of the discrete space is 1).
The key property of this measure that we use is that, if $U_\l \defined\set{b\in\cube^n \st b_\l = 1}$ is an ``upper facet'' of the hypercube, the total measure of $P_n$ in any $U_\l$ is only slightly less than $1/2$.

With this groundwork laid, we will have a density $P_n$ on $\RR^d \times \cube^n$ with almost half the mass in any upper facet.
We now iteratively pull this density back down again, producing densities $Q_k$ on $\RR^d \times \cube^k$, starting from $Q_n = P_n$ and finishing with $Q_0$ a density on $K_0$.
When considering the first step of the lifting procedure, we noted that mass being at a point $x$ in the upper facet meant that it could be sent to either $x-av_1$ or $x+av_1$.
This holds in general, and so, if almost half the mass of $Q_k$ is in the $k$th upper facet, then we can choose to pull this mass back either adding or subtracting $av_k$ to/from the position at each point.
This way, even if the mass in the lower facet is pulled back adversarially, we can still exert enough control that $\mu(Q_{k-1}) = \mu(Q_k) + \eps_k v_k$ for some $\eps_k\in\signs$.
This leaves $\mu(Q_0) = \sum_i \eps_i v_i$, as required.

\bigskip

Then, in \Cref{sec:strong-counter}, we prove \Cref{thm:false}, showing that the strong Koml\'os conjecture is false. 
This proof is inductive in nature, and we induct on a slightly stronger statement than the one in the theorem.
To be precise, we say that a sequence $v_1,\dotsc,v_n$ of vectors \emph{forces height $a$} if, for any starting vector $p$ and any sequence $\eps_1,\dotsc,\eps_n$ of signs, there is some $k$ such that $\norm[\big]{p + \sum_{i=1}^k \eps_i v_i}_\infty \geq a$; the inclusion of the starting vector $p$ allows us to more easily combine these sequences together in series.

We then construct such sequences: we first construct a sequence which can force several entries to be large simultaneously, and then construct a further sequence which gives us more control over which entries are large simultaneously. 
In particular, we would like to partition the set of coordinate directions into pairs so that the entries corresponding to some pair of coordinates are forced to be large simultaneously. 
Once we have such a sequence of vectors, we can then apply a rotation of the coordinate axes: each pair of coordinate directions encodes a plane, and we rotate each plane by $\pi/4$ so that a vector being large in both directions before the rotation corresponds to it being even larger in a single coordinate direction after the rotation.


\subsection{Statement on AI use}
The proofs presented here were all originally found with the assistance of AI by various groups of authors (as detailed in \Cref{subsec:ad}) and the present author used various AI tools to help understand, simplify, and break down the proofs.
However, the exposition presented here was written entirely by the author after understanding the material fully. 
AI was directly used in the production of this manuscript purely for proofreading, checking for correctness, and writing code for tikz diagrams.


\section{Proof of the Koml\'os Conjecture}
\label{sec:back}

We first provide a motivated and simplified version of the proof, cutting out some details and making some assumptions, before we proceed to the formal version, which is presented in \Cref{subsec:collapse-proof}.


\subsection{Proof motivation}
\label{subsec:collapse-motivation}

Let $K_0 = (-C,C)^d$ and recall that we write $\mu(P) = \int x \, P(x) \d x$ for the centre of mass or mean of a probability density.
Our aim is to construct a probability density $P$ on $K_0$ with centre of mass at $\mu(P) = \sum_i \eps_i v_i$ for some signs $\eps_i = \pm 1$.
We will do this in steps, producing distributions with mean at $\sum_{i=k}^n \eps_i v_i$ for each $k=n,n-1,\dotsc,1$ in turn; the reason for the reverse ordering will be explained in due course.
Performing just one step of this process is not too difficult; to give a trivial example, we could put all mass of $P$ near the centre of $K_0$, and then the first step causes no problems: we can just move all the mass by $+v_1$.
The difficulty is how to make many such steps in sequence. 

First of all, it will be helpful to have some control over just how we are going to go about distorting the probability density to move its mean. 
In particular, when we are trying to move the mass by $\pm v_i$, it makes sense to only move mass from points $x$ to points $x + \lambda v_i$ for some real $\lambda$.
In fact, we can go further, and fix some number $a > 0$ (think, for example, of $a=5$), and only move mass from $x$ to one of the points $x$, $x - av_i$, and $x + av_i$.
If a large enough proportion of the mass is sent to $x + av_i$, we will move the mean by $+v_i$.

\begin{figure}[htb]
    \centering
\begin{tikzpicture}[
  x=1cm,y=1cm,line width=.7pt,line cap=round,line join=miter,
  >=Stealth,font=\normalsize,
  outline/.style={draw=black,line width=.8pt,line cap=butt,line join=miter},
  point/.style={cross out,draw,minimum size=4pt,inner sep=0pt},
  guide/.style={densely dashed,line width=.5pt}]
  \draw[outline] (-1.4,-1.4) rectangle (1.4,1.4);
  \node[point] at (-.85,-.60) {};
  \node[anchor=west] at (-.73,-.60) {$x-av_1$};
  \node[point] at (.25,.50) {};
  \node[anchor=south] at (.25,.70) {$x+av_1$};
  \node at (0,-1.85) {$K_0$};
  \node at (10.85,-1.85) {$K_0$};

  \draw[->] (1.80,0) -- (3.05,0)
    node[midway,above=7pt] {Lift};

  \begin{scope}[xshift=5.25cm]
    \draw[outline] (-1.95,-1.95) -- (.85,-1.95) -- (.85,-.85)
      -- (1.95,-.85) -- (1.95,1.95) -- (-.85,1.95)
      -- (-.85,.85) -- (-1.95,.85) -- cycle;
    \draw[guide] (-.85,-.85) rectangle (.85,.85);
    \node[point] at (-.30,-.05) {};
    \node[anchor=west] at (-.18,-.05) {$x$};
    \node at (.25,-.50) {$L$};
    \node at (0,-2.45) {$(K_0-av_1)\cup(K_0+av_1)$};
  \end{scope}

  \draw[->] (7.65,0) -- (8.95,0)
    node[midway,above=7pt] {Pullback};
  \begin{scope}[xshift=10.75cm]
    \fill[pattern=north east lines,pattern color=black!65]
      (-1.4,-1.4) -- (.3,-1.4) -- (.3,-.3) -- (1.4,-.3)
      -- (1.4,1.4) -- (-.3,1.4) -- (-.3,.3) -- (-1.4,.3) -- cycle;
    \draw[guide] (-1.4,.3) -- (-.3,.3) -- (-.3,1.4);
    \draw[guide] (.3,-1.4) -- (.3,-.3) -- (1.4,-.3);
    \draw[outline] (-1.4,-1.4) rectangle (1.4,1.4);
  \end{scope}
\end{tikzpicture}
	\caption{A simple display of lifting and pulling back a density. If $P$ is a density on $K_0$, then the mass at $x\in L$ in the lifted density is $\frac{1}{2}\p{P(x-av_1) + P(x+av_1)}$. After pulling back, the mass in the shaded region might have changed, whereas the mass in the unshaded region will be the same as at the start. The idea is that, by performing this operation, we can move the mean of the distribution by $\pm v_1$.}
	\label{fig:simple-lift}
\end{figure}
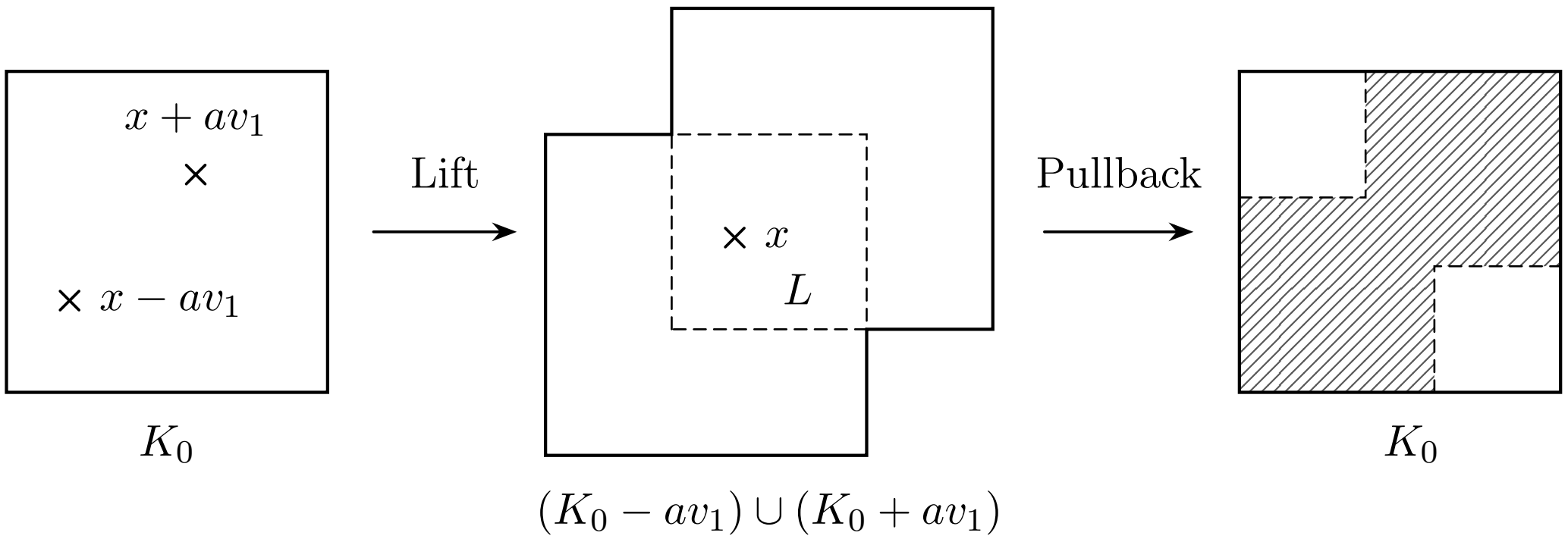

How exactly will we go about doing this?
There is a convenient framing of this way of moving mass around in terms of what we call ``lifts'' and ``collapses''.
We now discuss a toy example of this process.
The sets we define for now are not the final definitions; the argument will not work, and from this we will see how to change the definitions.

Let $K_1 = (K_0 - av_1) \union (K_0 + av_1)$.
We can naturally lift a probability density $P_0$ supported on $K_0$ (but defined on all of $\RR^d$, zero outside $K_0$) to $P_1$ on $K_1$ by setting $P_1(x) = \frac{1}{2}\p{P_0(x - av_1) + P_0(x + av_1)}$.
We can then pull this density back to a different density $P_0^*$ on $K_0$: consider the set $L = (K_0 - av_1) \inter (K_0 + av_1) \sseq K_1$.
Any mass of $P_1$ at a point $x\in L$ will have had contributions from both $x-av_1$ and $x+av_1$.
When pulling back to $P_0^*$, we will send all of this mass back to just one of these two points.
In particular, for all points in $L$, we have a choice of where to send the mass, and on $K_1 \setminus L$, we have no choice, and must send the mass to whichever of $x-av_1$ and $x+av_1$ is inside $K_0$.
This power of choice in $L$ is exactly what we will use to move the mean.
\Cref{fig:simple-lift} shows an example of this.

One thing that may be apparent is that the collapse operation described above is in some ways quite unpleasant: it can take a smooth distribution and make it extremely non-smooth (if we are not careful, it can even fail to be measurable!) and so it is hard to see how we could preserve any nice property of a distribution through this operation.
This brings us to the next key idea: we perform all the lift steps first, and then all the collapses. 

This is not how we will define the domains in the final proof, but for now iterate the process shown in \Cref{fig:simple-lift} and take $K_k = (K_{k-1} - av_k) \union (K_{k-1} + av_k)$ and $L_k = (K_{k-1} - av_k) \inter (K_{k-1} + av_k)$. 
We can also naturally lift a density $P_0$ on $K_0$ to $P_k$ on $K_k$ with $P_k(x) = \frac{1}{2}\p{P_{k-1}(x - av_k) + P_{k-1}(x + av_k)}$.
Note that the centre of mass of $P_n$ is the same as that of $P_0$, i.e.\ $\mu(P_n) = \mu(P_{n-1}) = \dotsb = \mu(P_0)$; it is only the collapse steps which move this centre of mass.
Let $Q_n = P_n$. We will now construct densities $Q_{n-1}, Q_{n-2},\dotsc,Q_0$ in sequence, with $Q_k$ supported on $K_k$, such that the centre of mass of $Q_{k-1}$ is equal to the centre of mass of $Q_{k}$ plus $\eps_k v_k$; the density $Q_{k-1}$ will be a distortion of $Q_k$.

The collapse steps distort the measures considerably.
However, assume for now that, when we get to pulling back $Q_k$ to $Q_{k-1}$, we have a significant amount of mass of $Q_k$ inside $L_k$.
This mass at a point $x\in L_k$ can be pulled back to either $x-av_k$ or $x+av_k$; as there is a large quantity of such mass, we have a lot of control over the mean of $Q_{k-1}$.
In particular, we need enough mass to overcome the fact that we have no control over where the mass on $K_k \setminus L_k$ will be sent: control over just less than half of the mass will suffice.
We will return to the problem of pulling back in due course.

The question, then, is: how can we force there to be a significant amount of mass in $L_k$?
Note that, if the distributions $P_0 - av_1$ and $P_0 + av_1$ are very similar (i.e.\ the total variation distance between them is small), then a lot of mass is in the overlap of the two distributions.
This overlap is precisely the set $L_1$. 
The distance between these two translates will be called the \emph{shift distance} of the measure, and we will find a way to make sure that the shift distance does not increase as we lift the distribution.
It is for this reason that we actually lift to a more complex set than that described above.

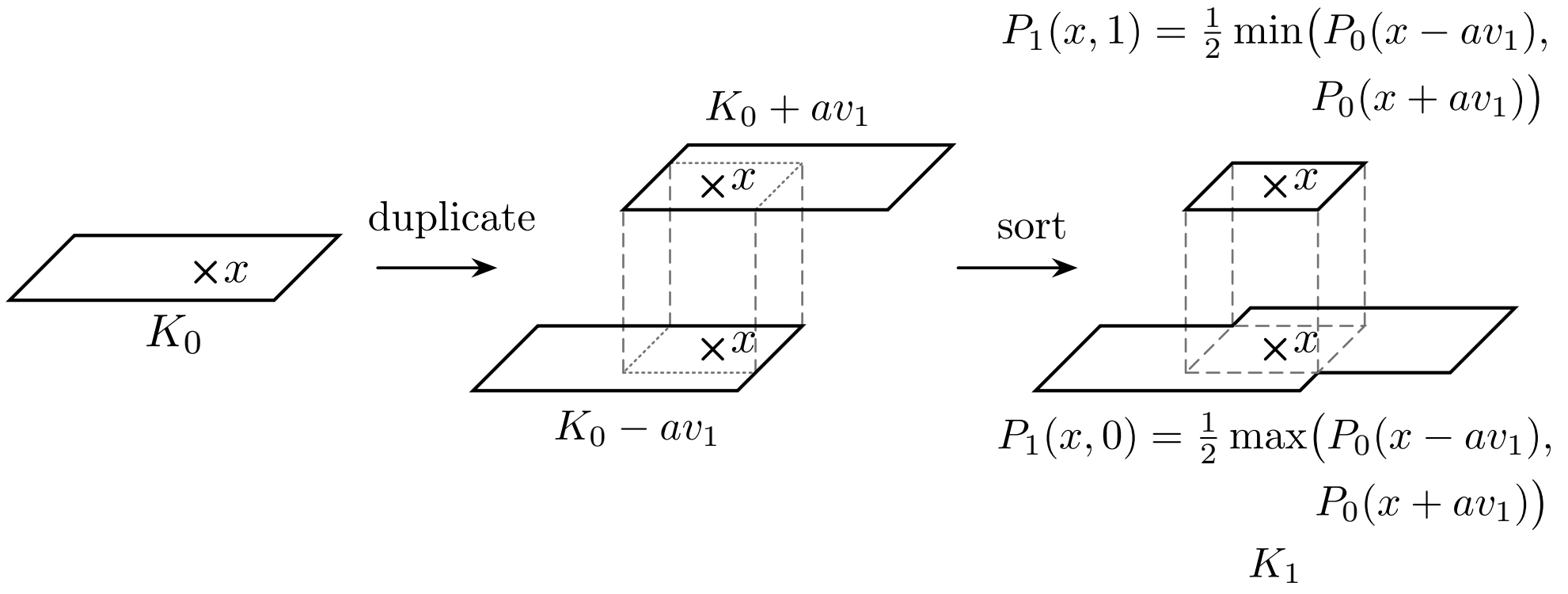
\begin{figure}[htb]
    \centering
\begin{tikzpicture}[
  x={(.85cm,0cm)},y={(.3005204cm,.3005204cm)},z={(0cm,.85cm)},
  line width=.7pt,line cap=round,line join=miter,
  >=Stealth,font=\normalsize,
  outline/.style={draw=black,line width=.8pt,line cap=butt,line join=miter},
  point/.style={cross out,draw,minimum size=4pt,inner sep=0pt},
  guide/.style={densely dashed,black!55,line width=.5pt},
  density/.style={font=\small,align=center,inner sep=2pt}]
  \begin{scope}[yshift=.68cm]
    \draw[outline] (-1.30,-.90,0) -- (1.30,-.90,0)
      -- (1.30,.90,0) -- (-1.30,.90,0) -- cycle;
    \node[point] at (.35,-.12,0) {};
    \node[anchor=west,inner sep=1pt] at (.47,-.12,0) {$x$};
    \node at (0,0,-.65) {$K_0$};
  \end{scope}

  \begin{scope}[xshift=4.5cm]
    \foreach \xx in {-.65,.65} {
      \foreach \yy in {-.65,.65} {
        \draw[guide] (\xx,\yy,0) -- (\xx,\yy,1.60);
      }
    }
    \foreach \zz in {0,1.60} {
      \draw[guide,densely dotted] (-.65,-.65,\zz) -- (.65,-.65,\zz)
        -- (.65,.65,\zz) -- (-.65,.65,\zz) -- cycle;
    }
    \draw[outline] (-1.95,-1.15,0) -- (.65,-1.15,0)
      -- (.65,.65,0) -- (-1.95,.65,0) -- cycle;
    \draw[outline] (-.65,-.65,1.60) -- (1.95,-.65,1.60)
      -- (1.95,1.15,1.60) -- (-.65,1.15,1.60) -- cycle;
    \node[point] at (0,0,0) {};
    \node[point] at (0,0,1.60) {};
    \node[anchor=west,inner sep=1pt] at (.12,0,.07) {$x$};
    \node[anchor=west,inner sep=1pt] at (.12,0,1.67) {$x$};
    \begin{scope}[x={(1cm,0cm)},y={(0cm,1cm)}]
      \node[density,anchor=south] at (.63,1.78)
        {$K_0+av_1$};
      \node[density,anchor=north] at (-.63,-.43)
        {$K_0-av_1$};
    \end{scope}
  \end{scope}

  \begin{scope}[xshift=9.2cm]
    \foreach \xx in {-.65,.65} {
      \foreach \yy in {-.65,.65} {
        \draw[guide] (\xx,\yy,0) -- (\xx,\yy,1.60);
      }
    }
    \draw[outline] (-1.95,-1.15,0) -- (.65,-1.15,0)
      -- (.65,-.65,0) -- (1.95,-.65,0) -- (1.95,1.15,0)
      -- (-.65,1.15,0) -- (-.65,.65,0) -- (-1.95,.65,0) -- cycle;
    \draw[guide] (-.65,-.65,0) -- (.65,-.65,0)
      -- (.65,.65,0) -- (-.65,.65,0) -- cycle;
    \draw[outline] (-.65,-.65,1.60) -- (.65,-.65,1.60)
      -- (.65,.65,1.60) -- (-.65,.65,1.60) -- cycle;
    \node[point] at (0,0,0) {};
    \node[point] at (0,0,1.60) {};
    \node[anchor=west,inner sep=1pt] at (.12,0,.07) {$x$};
    \node[anchor=west,inner sep=1pt] at (.12,0,1.67) {$x$};
    \begin{scope}[x={(1cm,0cm)},y={(0cm,1cm)}]
      \node[density,anchor=south] at (0,1.78)
        {$\begin{aligned}
           P_1(x,1)&=\tfrac12\min\bigl(P_0(x-av_1),\\
             &\hphantom{=\tfrac12\min\bigl(}P_0(x+av_1)\bigr)
         \end{aligned}$};
      \node[density,anchor=north] at (0,-.43)
        {$\begin{aligned}
           P_1(x,0)&=\tfrac12\max\bigl(P_0(x-av_1),\\
             &\hphantom{=\tfrac12\max\bigl(}P_0(x+av_1)\bigr)
         \end{aligned}$\\[3pt]$K_1$};
    \end{scope}
  \end{scope}

  \begin{scope}[x={(1cm,0cm)},y={(0cm,1cm)}]
    \draw[->] (1.70,.68) -- (2.70,.68)
      node[midway,above=3pt,xshift=1.2mm,font=\small] {duplicate};
    \draw[->] (6.55,.68) -- (7.55,.68)
      node[midway,above=3pt,xshift=1.2mm,font=\small] {sort};
  \end{scope}
\end{tikzpicture}
	\caption{A schematic diagram of how $K_1$ and $P_1$ are actually constructed. It is shown in a three dimensional projection to give the impression that the upper and lower parts are above and below each other (whereas they are strictly in different copies of $\RR^2$). The values of the density at points above $x$ are labelled.}
	\label{fig:schematic-p1}
\end{figure}

At this point it is convenient to start giving concrete definitions which match those used in the formal proofs: the shift distance $\Delta(P,v)$ is defined as
\begin{equation*}
	\Delta(P,v) \defined \dTV(P,P+v) = \frac{1}{2} \int \abs[\big]{P(x) - P(x - v)}\d x.
\end{equation*}
If $\Delta(P_0,2av_1) < \delta$, then $P_0 - av_1$ and $P_0 + av_1$ are very similar, and so $P_1$ has mass at least $1-\delta$ in $L_1$.
We need the lift process to not increase the shift distance, and we need to reserve almost half of the mass of $P_1$ in such a way that it will always be in $L_1$, even after further lift and collapse steps.
The key to doing this is the construction shown in \Cref{fig:schematic-p1}.
This is the actual construction that we use in the proof.

Taking two copies of $P_0$ in different spaces (each having its mass scaled down by a factor of $2$ so the combination remains normalised) and considering them as a single density does not change the shift distance.
This step is labelled ``duplicate'' in \Cref{fig:schematic-p1}. 

\begin{figure}[htb]
    \centering
\begin{tikzpicture}[
  x=1cm,y=1cm,line width=.7pt,line cap=round,line join=miter,
  font=\normalsize,
  outline/.style={draw=black,line width=.8pt,line cap=butt,line join=miter},
  point/.style={cross out,draw,minimum size=4pt,inner sep=0pt},
  density/.style={font=\small,fill=white,inner sep=6pt}]
  \draw[densely dashed,black!60,line width=.55pt]
    (0,-5.95) -- (0,4.70);
  \draw[densely dashed,black!60,line width=.55pt]
    (-6.65,0) -- (6.65,0);

  \begin{scope}[shift={(-3.45,2.65)}]
    \node[anchor=north west] at (-2.90,1.90) {$(1,0)$};
    \draw[outline] (-.95,-.40) -- (-.35,-.40) -- (-.35,-.90)
      -- (.95,-.90) -- (.95,.40) -- (.35,.40)
      -- (.35,.90) -- (-.95,.90) -- cycle;
    \node[point] at (0,0) {};
    \node[anchor=west,inner sep=1pt,font=\small] at (.06,0) {$x$};
    \node[density] at (0,-1.65)
      {$\begin{aligned}
         P_2(x,1,0)&=\tfrac12\max\bigl(P_1(x-av_2,1),\\
           &\hphantom{=\tfrac12\max\bigl(}P_1(x+av_2,1)\bigr)
       \end{aligned}$};
  \end{scope}

  \begin{scope}[shift={(3.45,2.65)}]
    \node[anchor=north east] at (2.90,1.90) {$(1,1)$};
    \draw[outline] (-.35,-.40) rectangle (.35,.40);
    \node[point] at (0,0) {};
    \node[anchor=west,inner sep=1pt,font=\small] at (.06,0) {$x$};
    \node[density] at (0,-1.65)
      {$\begin{aligned}
         P_2(x,1,1)&=\tfrac12\min\bigl(P_1(x-av_2,1),\\
           &\hphantom{=\tfrac12\min\bigl(}P_1(x+av_2,1)\bigr)
       \end{aligned}$};
  \end{scope}

  \begin{scope}[shift={(-3.45,-2.90)}]
    \node[anchor=north west] at (-2.90,1.90) {$(0,0)$};
    \draw[outline] (-2.25,-.90) -- (-1.65,-.90) -- (-1.65,-1.40)
      -- (.95,-1.40) -- (.95,-.90) -- (2.25,-.90)
      -- (2.25,.90) -- (1.65,.90) -- (1.65,1.40)
      -- (-.95,1.40) -- (-.95,.90) -- (-2.25,.90) -- cycle;
    \node[density] at (0,-2.15)
      {$\begin{aligned}
         P_2(x,0,0)&=\tfrac12\max\bigl(P_1(x-av_2,0),\\
           &\hphantom{=\tfrac12\max\bigl(}P_1(x+av_2,0)\bigr)
       \end{aligned}$};
    \node[point] at (0,0) {};
    \node[anchor=west,inner sep=1pt,font=\small] at (.06,0) {$x$};
  \end{scope}

  \begin{scope}[shift={(3.45,-2.90)}]
    \node[anchor=north east] at (2.90,1.90) {$(0,1)$};
    \draw[outline] (-1.65,-.90) -- (.35,-.90) -- (.35,-.40)
      -- (1.65,-.40) -- (1.65,.90) -- (-.35,.90)
      -- (-.35,.40) -- (-1.65,.40) -- cycle;
    \node[point] at (0,0) {};
    \node[anchor=west,inner sep=1pt,font=\small] at (.06,0) {$x$};
    \node[density] at (0,-2.15)
      {$\begin{aligned}
         P_2(x,0,1)&=\tfrac12\min\bigl(P_1(x-av_2,0),\\
           &\hphantom{=\tfrac12\min\bigl(}P_1(x+av_2,0)\bigr)
       \end{aligned}$};
  \end{scope}
  \node[fill=white,inner sep=4pt] at (0,-6.30) {$K_2$};
\end{tikzpicture}
	\caption{A diagram of how $K_2$ and $P_2$ are constructed from $K_1$ and $P_1$. We have switched to viewing the four separate parts of this distribution as in four different copies of $\RR^2$.}
	\label{fig:schematic-p2}
\end{figure}
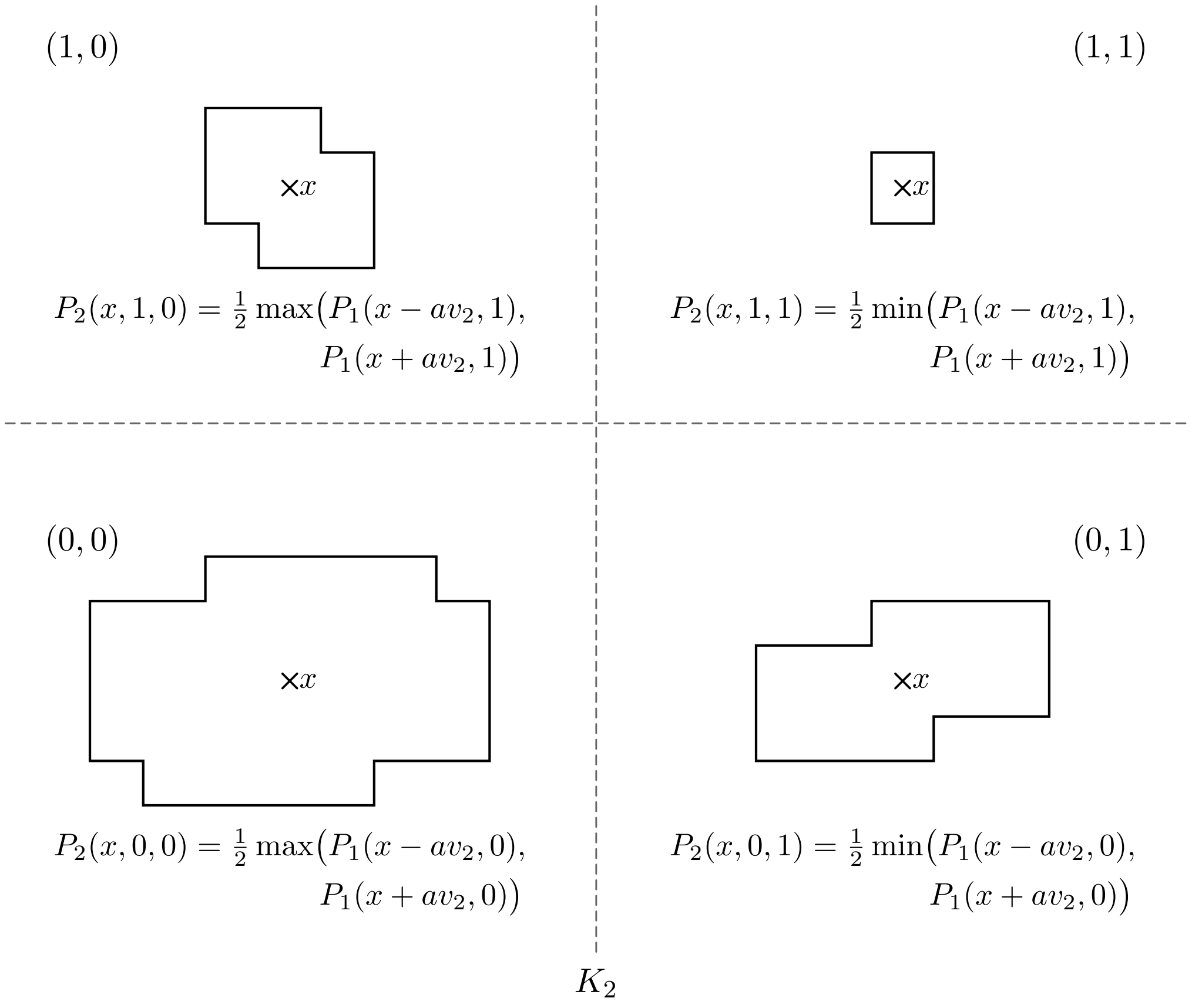

The ``sort'' step does not increase shift distance (but may decrease it, which only helps us). 
The intuition behind this is that, after sorting, computing the shift distance compares large values with large values and small values with small values, which can only decrease the differences involved.
If the shift distance is sufficiently small, then the distributions on the two parts of $K_1$ are in fact very similar.
As it is not immediately clear how to generalise this construction to further steps, a schematic of $K_2$ is shown in \Cref{fig:schematic-p2}: there are now four parts of the density in different copies of $\RR^d$.

In general we can define $K_k \sseq \RR^d\times\cube^k$ and $P_k$ a density supported on $K_k$ (and defined to be zero outside this set) by iterating: the process shown in \Cref{fig:schematic-p1} is applied to the parts of $K_{k-1}$ and $P_{k-1}$ in each copy of $\RR^d$ to produce $K_k$ and $P_k$.

We remark here that, while we have labelled the elements of the discrete space $\cube^k$ with $0$ and $1$, we are thinking of these as abstract labels rather than numbers: we could happily run the whole proof with space called e.g.\ $\set{\bot,\top}^k$.
This part of the space really is just a discrete hypercube.
In particular, it will not interact with arithmetic operations in any way: we will abuse notation and write e.g.\ $x + v$ where $x\in\RR^d\times\cube^k$ and $v\in\RR^d$, and this will mean to simply add $v$ to the $\RR^d$ part of $x$ (if $x = (z, b)$ with $z\in\RR^d$, then $x+v = (z+v, b)$).

It was promised earlier in the outline that the centre of mass of the distributions would satisfy $\mu(P_0) = \mu(P_1) = \dotsb = \mu(P_n)$.
How does this track with the discrete coordinates being added in?
Following the above discussion about the discrete coordinates not being thought of as numbers, we ignore them for the purposes of computing the mean.
Define $\pi\from\RR^d\times\cube^k \to \RR^d$ to be the projection onto $\RR^d$; we can then define the centre of mass of a density $P_k$ on $\RR^d\times\cube^k$ to be
\begin{equation*}
	\mu(P_k) \defined \int_{\RR^d\times\cube^k} \pi(x) \, P_k(x) \d x.
\end{equation*}
With this definition in hand, it should then be clear from observing \Cref{fig:schematic-p1} that $\mu(P_1) = \mu(P_0)$: the ``duplicate'' step is symmetric and the ``sort'' step only changes discrete coordinates.
Similarly, we find that no lift step changes the centre of mass of a density, and so $\mu(P_n) = \mu(P_0) = 0$.

In the toy example we considered earlier, we said that we needed almost half the mass to be in the overlap set $L$. What does that mean here?
Note that, once a discrete coordinate is set for some mass, it does not change: the mass in the ``upper'' part of \Cref{fig:schematic-p1} will be moved around as we perform further lifts, but its first discrete coordinate in $\cube^k$ will always remain $1$.
With this in mind, let the \lth\ upper facet of the hypercube $\cube^k$ be the set of those points $b\in\cube^k$ with $b_\l = 1$.

The total mass of $P_k$ in the \lth\ upper facet of $\cube^k$ is equal to the total mass of $P_\l$ in the \lth\ upper facet of $\cube^\l$: these discrete coordinates are not changed by later operations. 
That means, because the shift distance of $P_k$ in any direction is at most the shift distance of $P_0$, which is small, $P_k - av_{k+1}$ and $P_k + av_{k+1}$ are very similar distributions.
This means that their minimum and maximum are very similar (think of the extreme case: if, somehow, the two shifted densities were equal, then the min and max would be equal), and so almost half the mass of $P_{k+1}$ is contained in the ``minimum'' part---that is, the upper facet.

We therefore find that every upper facet of $\cube^n$ contains almost half the mass of $P_n$ (note in particular that different upper facets have considerable intersections with each other).
This is the key property that we will use in pulling back the densities.

Let $Q_n = P_n$; we iteratively define $Q_{k-1}$ supported on $K_{k-1}$ in terms of $Q_k$.
When working with $Q_k$, we would like to still know that there is lots of mass in the \kth\ upper facet.
It is therefore important that our collapses do not change any discrete coordinates besides the one they are getting rid of.
The collapses follow essentially the same intuition as in the toy example at the start of this section: mass at a point $(x,b) \in \RR^d\times\cube^k$ in the \kth\ upper facet can be pulled back to either $p_- = (x - av_k, b')$ or $p_+ = (x + av_k, b')$, where $b'\in\cube^{k-1}$ is $b$ with the final coordinate removed.
Note that both $p_+$ and $p_-$ are in $K_{k-1}$.

However we pull the mass back, we certainly have that $\mu(Q_{k-1}) = \mu(Q_k) + \lambda v_k$ for some real $\lambda \in [-a, a]$.
Our goal is to hit $\lambda = \pm 1$.
Assume that the mass outside the \kth\ upper facet (i.e.\ in the \kth\ lower facet) is pulled back by an adversary.
By pulling back just over half of the mass, the adversary can change $\lambda$ by slightly more than $a/2$.
But the distance from such a point to the closer of $-1$ and $+1$ is significantly less than $a/2$, which is an amount of control that we can exert.
Therefore we can arrange for the collapse to satisfy $\mu(Q_{k-1}) = \mu(Q_k) + \eps_k v_k$ for some $\eps_k = \pm 1$, as required.

Iterating this procedure, we will eventually produce a density $Q_0$ back on $K_0$ with mean equal to $\sum_i \eps_i v_i$, proving the Koml\'os conjecture.


\subsection{Formal proof}
\label{subsec:collapse-proof}

We now present the formal details of the proof of \Cref{thm:main}, referring back to \Cref{subsec:collapse-motivation} to motivate the ideas involved.
We first give formal definitions of the terms we will use in the proof, and then move on to a series of lemmas.
Throughout the entire proof, the values of $d$ and $n$ are fixed.

\begin{definition}
	A set $X$ is called \emph{nice} if we have $X\sseq \RR^d \times \cube^k$ for some integer $k$ such that, for all $b\in \cube^k$, the fibre
	$\set[\big]{x\in\RR^d \st (x,b) \in X}$ is an open bounded set.

	We will abuse notation by writing e.g.\ $x + v$ when $x \in X$ and $v \in \RR^d$: if $x = (u, b)$ with $u \in \RR^d$ and $b \in \cube^k$, then $x + v$ is shorthand for $(u + v, b)$.

	We write $\pi\from X \to \RR^d$ and $\sigma \from X \to \cube^k$ for the projection maps
	\begin{equation*}
		\pi(x, b) \defined x \quad \text{and} \quad \sigma(x, b) \defined b \qquad \text{ for all } x \in \RR^d \text{ and } b \in \cube^k.
	\end{equation*}
	We will again abuse notation and not distinguish between such projection maps with different domains.

	Given a probability density $P$ on a nice set $X$, we write $\mu(P)$ for the centre of mass of $P$ after projecting by $\pi$. That is,
	\begin{equation*}
		\mu(P) \defined \int_X \pi(z) P(z) \d z = \sum_{b\in\cube^k} \int_{\RR^d} x \, P(x, b) \d x,
	\end{equation*}
	where we assume that $P$ takes value $0$ everywhere outside of $X$.
	
	We refer to a pair $(X, P)$ as \emph{nice} if it consists of a nice set $X$ and a probability density on $X$.
\end{definition}

We now recall the formal definition of the shift distance of a measure in a direction.

\begin{definition}
	Given two probability densities $P$ and $Q$ on a set $X \sseq \RR^d\times \cube^k$, the \emph{total variation distance} between $P$ and $Q$ is
	\begin{equation*}
		\dTV(P,Q) \defined \frac{1}{2} \int_{\RR^d\times \cube^k} \abs[\big]{P(x) - Q(x)} \d x.
	\end{equation*}
	We define the \emph{shifted} function $P+u$ for $u\in\RR^d$ by $(P+u)(x) \defined P(x - u)$.
	We can then define the \emph{shift distance} $\Delta(P,u)$ as follows.
	\begin{equation*}
		\Delta(P,u) \defined \dTV(P, P+u) = \frac{1}{2}\int_{\RR^d\times \cube^k} \abs[\big]{P(x) - P(x - u)} \d x.
	\end{equation*}
\end{definition}

Given a probability density $Q$ on a set $X$ and a subset $Y \sseq X$, we write $Q[Y]$ for the total density of $Y$, that is
\begin{equation*}
	Q[Y] \defined \int_Y Q(y) \d y.
\end{equation*}

With these definitions in hand, we can now proceed to the first main section of this proof: lifting.


\subsubsection{Lifts}
\label{subsubsec:lift}

\begin{definition}[$v$-lifts]
	Given a nice set $X \sseq \RR^d \times \cube^k$, a probability density $P$ supported on $X$, and a vector $v\in\RR^d$, the \emph{$v$-lift} of $(X,P)$ is a pair $(X',P')$ where $X' \sseq \RR^d \times \cube^{k+1}$ and $P'$ is a probability density supported on $X'$. 
	These objects are defined as follows.
	\begin{equation*}
		X' \defined \b[\big]{((X - v) \union (X + v)) \times \set{0}} \union \b[\big]{((X - v) \inter (X + v)) \times \set{1}}.
	\end{equation*}
	The \emph{$v$-lifted distribution} $P'$ is defined, for all $x \in \RR^d \times \cube^k$ by
	\begin{align*}
		P'(x,0) &\defined  \tfrac{1}{2} \max\set[\big]{P(x - v), \, P(x + v)}, \\
		P'(x,1) &\defined  \tfrac{1}{2} \min\set[\big]{P(x - v), \, P(x + v)}.
	\end{align*}
	The \emph{$v$-lift} of $(X,P)$ is the pair $(X',P')$ consisting of the $v$-lifts of the two components.
\end{definition}

See \Cref{fig:schematic-p1,fig:schematic-p2} for diagrams of a lift.
We note before continuing that if $(X',P')$ is the $v$-lift of a nice pair $(X,P)$, then $(X',P')$ is also a nice pair.
It is immediate that $X'$ is a nice set, and the fact that $P'$ is a probability density function follows from the simple identity
\begin{equation}
	\label{eq:simple-sum}
	P'(x, 0) + P'(x, 1) = \tfrac{1}{2} \p[\big]{P(x - v) + P (x + v)}.
\end{equation}

Another basic property of the $v$-lift which follows from \eqref{eq:simple-sum} is that it does not change the mean of (the projection to $\RR^d$ of) the density, as follows.

\begin{lemma}
	\label{lem:lift-preserves-mean}
	If $(X,P)$ is a nice pair with $v$-lift $(X',P')$, then $\mu(P') = \mu(P)$.
\end{lemma}

This fact should be clear upon considering \Cref{fig:schematic-p1}: the ``duplicate'' step in that diagram takes two copies of $P$ (each scaled down by a factor of $2$) and moves them in opposite directions, which does not change the mean. The ``sort'' step then only changes the discrete coordinates, and so again has no effect on the mean.
The following proof formalises this intuition.

\begin{proof}[Proof of \Cref{lem:lift-preserves-mean}]
	We may compute
	\begin{align*}
		\mu(P') 
		&= \int_{z'\in X'} P'(z') \pi(z') \d z' \\
		&= \int_{z\in \RR^d\times \cube^k} \b[\big]{P'(z,0) + P'(z,1)} \pi(z) \d z \\
		&= \int_{z\in \RR^d\times \cube^k} \frac{1}{2}\b[\big]{P(z - v) + P(z + v)}\pi(z) \d z \\
		&= \int_{z\in \RR^d\times \cube^k} \frac{1}{2} P(z) \b[\big]{\pi(z + v) + \pi(z - v)} \d z \\
		&= \int_{z\in \RR^d\times \cube^k} P(z) \pi(z) \d z \\
		&= \mu(P).
	\end{align*}
	In the third line of the above, we used equality \eqref{eq:simple-sum}, and in the fifth, we used the fact that
	\begin{equation*}
		\pi(z + v) + \pi(z - v) = \pi(z) + v + \pi(z) - v = 2\pi(z).
	\end{equation*}
	This proves \Cref{lem:lift-preserves-mean}.
\end{proof}

Another simple property of $v$-lifts we prove is that the $v$-lift of a density function preserves mass in fibres over the discrete part of the set $X$.

\begin{lemma}
	\label{lem:lift-mass-discrete-fibre}
	Let $(X,P)$ be a nice pair with $X \sseq \RR^d \times \cube^k$, and let $(X',P')$ be the $v$-lift.
	Then, for all $b\in\cube^k$, we have
	\begin{equation*}
		P'\b[\big]{\sigma^{-1}(b,0)} + P'\b[\big]{\sigma^{-1}(b,1)} = P\b[\big]{\sigma^{-1}(b)}.
	\end{equation*}
\end{lemma}

The intuition behind this lemma is that the lift operation does not change any existing discrete coordinates: mass at $(x,b)\in \RR^d\times\cube^k$ is moved to points of the form $(x',b,b_{k+1})$ for some $x'\in\RR^d$ and $b_{k+1}\in\cube$.

\begin{proof}[Proof of \Cref{lem:lift-mass-discrete-fibre}]
	We may compute
	\begin{align*}
		P'\b[\big]{\sigma^{-1}(b,0)} &+ P'\b[\big]{\sigma^{-1}(b,1)} \\
		&= \int_{x \in \RR^d} \p[\big]{P'(x,b,0) + P'(x,b,1)} \d x \\
		&= \int_{x \in \RR^d} \frac{1}{2}\p[\Big]{\min\set[\big]{P(x-v,b), P(x+v,b)} + \max\set[\big]{P(x-v,b), P(x+v,b)}} \d x \\
		&= \int_{x \in \RR^d} \frac{1}{2}\p[\big]{P(x-v,b) + P(x+v,b)} \d x \\
		&= \int_{x \in \RR^d} P(x,b) \d x \\
		&= P\b[\big]{\sigma^{-1}(b)},
	\end{align*}
	as required.
\end{proof}

We now establish the key link between lifts and shift distance. 
We first show that the lifting operation does not increase shift distance.

\begin{lemma}
	\label{lem:shift-distance-decreasing}
	If $X$ is a nice set and $P$ is a probability density on $X$ with $v$-lifts $X'$ and $P'$ respectively (for some $v\in \RR^d$), then, for all $u \in \RR^d$, we have
	\begin{equation*}
		\Delta(P', u) \leq \Delta(P, u).
	\end{equation*}
\end{lemma}

As discussed in \Cref{subsec:collapse-motivation}, the intuition for this fact can be seen in \Cref{fig:schematic-p1}: the ``duplicate'' step does not change the shift distance, and the ``sort'' step leads to us comparing large values with large values and small values with small values, which can only decrease the distance.
The following proof captures this intuition formally.

\begin{proof}[Proof of \Cref{lem:shift-distance-decreasing}]
	Note that, due to \Cref{lem:lift-mass-discrete-fibre}, it suffices to prove this result for a single discrete fibre; we may thus assume that $X \sseq \RR^d$.
	We claim that
	\begin{align*}
		4\Delta(P',u) 
		&= 2\int_{x\in \RR^d} \b[\Big]{\abs[\big]{P'(x,0) - P'(x-u, 0)} + \abs[\big]{P'(x,1) - P'(x-u,1)}} \d x \\
		&\leq \int_{x\in \RR^d} \b[\Big]{\abs[\big]{P(x+v) - P(x+v-u)} + \abs[\big]{P(x-v) - P(x-v-u)}} \d x \\
		&= 2\int_{x\in \RR^d} \abs[\big]{P(x-u) - P(x)} \d x \\
		&= 4\Delta(P,u).
	\end{align*}
	The only step in the above calculation which is not immediate is the inequality.
	We in fact claim that this inequality follows pointwise.
	For ease of notation, write
	\begin{equation*}
		r_- = P(x - v),\; r_+ = P(x + v),\; s_- = P(x-v-u), \text{ and } s_+ = P(x+v-u).
	\end{equation*}
	We must therefore prove that
	\begin{equation}
		\label{eq:sorting-inequality}
		\abs[\big]{\max\set{r_-,r_+} - \max\set{s_-, s_+}} + \abs[\big]{\min\set{r_-,r_+} - \min\set{s_-, s_+}} \leq \abs{r_+ - s_+} + \abs{r_- - s_-},
	\end{equation}
	but this follows from a simple case analysis, which we nevertheless detail.
	If $r_+ \geq r_-$ and $s_+ \geq s_-$, then the two sides of \eqref{eq:sorting-inequality} are equal and we are done; likewise if $r_- \geq r_+$ and $s_- \geq s_+$.
	Assume therefore that $r_+ > r_-$ and $s_- > s_+$; the other case is similar.
	We may assume without loss of generality that $r_+ \geq s_-$.
	There are three cases: $r_+ > r_- \geq s_- > s_+$, $r_+ \geq s_- \geq r_- \geq s_+$, and $r_+ \geq s_- > s_+ \geq r_-$.

	In the first case, equality holds in \eqref{eq:sorting-inequality}.
	This inequality will be strict in the second case if $s_- > r_-$, and it always holds strictly in the third case.
	In particular, inequality \eqref{eq:sorting-inequality} always holds.

	This completes the proof of \Cref{lem:shift-distance-decreasing}.
\end{proof}

The next lemma will capture the intuition that, if $\Delta(P,2v)$ is small, then the distributions $P - v$ and $P+v$ are very similar, and so $\min\set{P - v, P + v}$ should on average not be much smaller than $\max\set{P - v, P + v}$.
In the language of lifts, there is a lot of mass in the ``upper'' part of the set.
To describe this all formally, we will introduce some notation for describing this ``upper'' part.
We let $\upstairs{k}{\l}$ and $\downstairs{k}{\l}$ be the sets of points in $\RR^d \times\cube^k$ having their $\l$th discrete coordinates equal to $1$ and $0$ respectively; more precisely, we make the following definition.

\begin{definition}
	Let $\upstairs{k}{\l} \sseq \RR^d \times\cube^k$ be the $\l$th \emph{upper facet} of space:
	\begin{equation*}
		\upstairs{k}{\l} \defined \set[\big]{(x_1,\dotsc,x_d,b_1,\dotsc,b_k) \in \RR^d\times\cube^k \st b_\l = 1}.
	\end{equation*}
	Similarly, define $\downstairs{k}{\l} \defined (\RR^d \times \cube^k) \setminus \upstairs{k}{\l}$ to be the $\l$th \emph{lower facet} of space.
\end{definition}

With this notation in hand, we can now formally state the following lemma.

\begin{lemma}
	\label{lem:lift-mass-upstairs}
	If $(X,P)$ is a nice pair on $\RR^d \times \cube^k$ with $v$-lift $(X',P')$ for some $v\in \RR^d$, then
	\begin{equation}
		\label{eq:lift-mass-upstairs}
		P'\b[\big]{\upstairs{k+1}{k+1}} = \frac{1 - \Delta(P, 2v)}{2}.
	\end{equation}
	Moreover, for all $\l \leq k$, we have
	\begin{equation}
		\label{eq:lift-mass-preserved}
		P'\b[\big]{\upstairs{k+1}{\l}} = P\b[\big]{\upstairs{k}{\l}}.
	\end{equation}
\end{lemma}

\begin{proof}
	Note that, as $v$-lifts preserve existing discrete coordinates, equality \eqref{eq:lift-mass-preserved} follows immediately.
	Moreover, in proving \eqref{eq:lift-mass-upstairs}, we may work separately in each discrete fibre and sum the results.
	We therefore assume that $X\sseq \RR^d$ (i.e.\ $k = 0$).
	We use the simple identity $2\min\set{a,b} = a + b - \abs{a-b}$, which holds for all $a,b\in \RR$, and compute
	\begin{align*}
		4P'\b[\big]{\upstairs{k+1}{k+1}} 
		&=  2 \int_{x\in \RR^d} \min\set[\big]{P(x+v), P(x-v)} \d x \\
		&= \int_{x\in \RR^d} \b[\Big]{P(x+v) + P(x-v) - \abs[\big]{P(x+v) - P(x-v)}} \d x \\
		&= 2\int_{x\in \RR^d}P(x) \d x - \int_{x\in \RR^d} \abs[\big]{P(x) - P(x - 2v)} \d x \\
		&= 2 - 2\Delta(P, 2v),
	\end{align*}
	as required.
\end{proof}

With \Cref{lem:shift-distance-decreasing,lem:lift-mass-upstairs} in hand, we now have sufficient control over the lifting part of the construction.
We now move on to the collapse half of the proof.


\subsubsection{Collapse}
\label{subsubsec:pull}

The first step is to formally define what a collapse actually is.

\begin{definition}
	Given a nice set $X\sseq \RR^d \times \cube^k$ and the $v$-lift $X'$ for some vector $v\in\RR^d$, a measurable function $f\from X' \to X$ is a \emph{$v$-collapse} function if, for all $(x,b_{k+1}) \in X'$ (where $x \in \RR^d \times \cube^k$ and $b_{k+1} \in \cube$), we have
	\begin{equation}
		\label{eq:collapse-definition}
		f(x,b_{k+1}) = x \pm v \in X.
	\end{equation} 
	Note that, by the definition of a $v$-lift, it is always possible to ensure that $f(x,b_{k+1}) \in X$.

	Then, given a probability density $P'$ on $X'$ ($P'$ not necessarily arising as a lift) and a $v$-collapse function $f$, the \emph{$f$-collapse} of $P'$ is the probability density $P^*$ on $X$ satisfying, for all $z\in X$,
	\begin{equation*}
		P^*(z) \defined \sum_{z'\in f^{-1}(z)} P'(z'),	
	\end{equation*}
	where we may note that $f^{-1}(z)$ contains at most four elements.
\end{definition}

Note that the value of $b_{k+1}$ is not used in the definition \eqref{eq:collapse-definition} of a $v$-collapse.
However, if $(x,1) \in X'$, then we know from the definition of the $v$-lift that both $x - v$ and $x + v$ are in $X$.
This is the power we will use: if there is a lot of mass in the upper facet, then we have control over whether this mass is collapsed to $x-v$ or $x+v$, and so have a lot of control over the mean of the distribution resulting from the collapse.

We first prove a simple property of collapse functions: the $v$-collapse does not change the mass in any fibers above $\cube^k$.

\begin{lemma}
	\label{lem:collapse-preserves-mass}
	Let $X\sseq \RR^d \times \cube^k$ be a nice set, and let $X'$ be the $v$-lift for some $v\in \RR^d$.
	Let $P'$ be an arbitrary measurable probability density on $X'$, and let $P^*$ be the $f$-collapse of $P'$, where $f$ is an arbitrary $v$-collapse of $X'$.
	Then, for any $b\in\cube^k$, we have
	\begin{equation*}
		P^*\b[\big]{\sigma^{-1}(b)} = P'\b[\big]{\sigma^{-1}(b,0)} + P'\b[\big]{\sigma^{-1}(b,1)}.
	\end{equation*}
\end{lemma}

\begin{proof}
	The key observation here is that the $v$-collapse function $f$ preserves the discrete component (except for the coordinate which is being collapsed):
	\begin{equation*}
		f^{-1}\b[\big]{\sigma^{-1}(b)} = \sigma^{-1}(b,0) \union \sigma^{-1}(b,1),
	\end{equation*}
	and there is no multiplicity.
	Therefore, we may compute
	\begin{align*}
		P^*\b[\big]{\sigma^{-1}(b)} 
		&= \int_{x\in \RR^d} P^*(x,b) \d x \\
		&= \int_{x\in \RR^d} \sum_{z\in f^{-1}(x,b)} P'(z) \d x \\
		&= \int_{z\in f^{-1}[\sigma^{-1}(b)]} P'(z) \d z \\
		&= \int_{x\in \RR^d} \p[\big]{P'(x,b,0) + P'(x,b,1)} \d x \\
		&= P'\b[\big]{\sigma^{-1}(b,0)} + P'\b[\big]{\sigma^{-1}(b,1)},
	\end{align*}
	as required.
\end{proof}

If one is concerned about measurability issues in the above computations, these are conceptually easy to resolve (though the notation would suffer): we can let $Q_n$ be a discretisation on a suitably fine mesh of $P_n$, and then all the integrals in the above proof become sums.

We can now state the lemma which gives us precise control over the mean of a collapse.
Recall that $\upstairs{k}{\l}$ and $\downstairs{k}{\l}$ are the sets of points in $\RR^d \times\cube^k$ having their $\l$th discrete coordinates equal to $1$ and $0$ respectively.

\begin{lemma}
	\label{lem:step}
	Let $X \sseq \RR^d\times\cube^k$ be a nice set and let $X'$ be its $v$-lift for some non-zero $v\in \RR^d$. 
	Write $U$ and $L$ for $\upstairs{k+1}{k+1}$ and $\downstairs{k+1}{k+1}$, the final upper and lower facets, respectively.
	Let $P'$ be an arbitrary probability density on $X'$ and assume that
	\begin{equation}
		\label{eq:collapse-mass-upstairs}
		P'\b{U} \geq \frac{1}{2} \max\set[\Big]{1 - \frac{1}{\norm{v}_2}, \; \frac{2}{\norm{v}_2}}.
	\end{equation}
	Then, for any measurable partial $v$-collapse function $f_0$ defined on $X'\inter L$, there is a total $v$-collapse function $f \from X' \to X$ extending $f_0$ such that, if $P^*$ is the $f$-collapse of $P'$, then
	\begin{equation}
		\label{eq:collapse-goal}
		\mu(P^*) = \mu(P') + \frac{\eps v}{\norm{v}_2}
	\end{equation}
	for some $\eps\in\set{-1, +1}$.
\end{lemma}

\begin{proof}
	Note first that, for any $v$-collapse function $f\from X' \to X$, if $P^*$ is the $f$-collapse of $P'$, then we have
	\begin{equation*}
		\mu(P^*) = \mu(P') + \lambda v \quad \text{ for some } \lambda \in [-1,1].
	\end{equation*}
	Note that we can choose where all mass in the upper facet $U$ is sent, sending the mass at each point $(x,i)\in U$ (for $i\in\cube$) to either $x+v$ or $x-v$.
	This means that there is a closed interval $I$ of width $2P'[U]$ with $I \sseq [-1,1]$ such that, for all $\lambda \in I$, there is a $v$-collapse function extending $f_0$ which moves the mean by $\lambda v$.
	To produce a $v$-collapse function satisfying \eqref{eq:collapse-goal}, we need $\lambda = \pm 1/\norm{v}_2$.
	These two values chop the interval $[-1,1]$ into three parts, two of length $1 - 1/\norm{v}_2$ and one of length $2/\norm{v}_2$.
	Thus, if $\abs{I}$ is at least these two values, wherever it is placed inside $[-1,1]$, it must hit one of them.

	The condition \eqref{eq:collapse-mass-upstairs} is precisely as required for $\abs{I} = 2 P'[U]$ to be large enough, and so the conclusion holds.
\end{proof}

This concludes the results we require about collapses; we can now proceed to show that the inductive construction can actually get off the ground.


\subsubsection{Constructing the starting probability density}
\label{subsubsec:construction}

The goal of this section is to construct a probability density on $K_0 = (-C,C)^d$ with small shift distance when shifted by any unit vector. 
It is not obvious that such a distribution should exist, and we make a couple of remarks on this topic before continuing to explain why such a density should exist.

\begin{remark}
	The uniform density on $K_0$ does not have small shift distance: fix $C>0$ and let $v = (1/\sqrt{d},\dotsc,1/\sqrt{d})$.
	Then
	\begin{equation*}
		\frac{\abs{K_0 \inter (K_0 + v)}}{\abs{K_0}} \to 0 \quad\text{ as }\quad d\to\infty.
	\end{equation*}
	In particular, the shift distance of the uniform density on $K_0$ in direction $v$ tends to 1 as $d$ tends to infinity.
\end{remark}

\begin{remark}
	This section is the only place in the proof in which we use the fact that the Koml\'os conjecture concerns the infinity norm; every other step requires only the existence of a suitable probability density.
	In particular, this means that one cannot find a probability density with small shift distance in every direction within a constant multiple (as $d\to\infty$) of the unit $L_p$ ball for any $p < \infty$.

	Indeed, the Koml\'os conjecture in the $L_p$ ball is easily seen to be false by considering the sequence of vectors $e_1,\dotsc,e_d$, the standard basis of $\RR^d$.
	In this case, all signed sums have $p$-norm equal to $d^{1/p}$.
\end{remark}

We essentially give two different proofs that there is a function with small shift distance in every direction on the hypercube. 
The first, which we outline informally, will show that many product functions have the required property and is, in the author's opinion, more intuitive than the version of the proof that we formalise (in which setting the second proof is superior).

We will take $P(x) = \prod_{i=1}^d g(x_i)$ for some non-negative symmetric function $g\from\RR\to\RR$ supported on $(-C,C)$ with $\int_\RR g(t) \d t = 1$.
Here we use the fact that we are working on the hypercube, as different coordinates behave independently.
In particular, picking $X\in(-C,C)^d$ randomly according to the density $P$ is equivalent to independently picking coordinates $X_1,\dotsc,X_d$ each according to $g$.

Pick $X$ randomly according to $P$ and consider the quantity $\partial_u P(X) / P(X)$, where $\partial_u = \sum_{i=1}^d u_i (\partial/\partial x_i)$ is the directional derivative in direction $u$.
If $\EE\abs[\big]{\partial_u P(X) / P(X)}$ is small, then we expect that upon taking a step of $\delta u$ from $X$ (for $\delta$ small), the relative value of $P$ does not change much.
This is exactly saying that $P$ and $P+\delta u$ are similar, i.e.\ $\Delta(P, \delta u)$ is small.
We can compute

\begin{equation*}
	\frac{\partial_u P(X)}{P(X)} = \frac{\sum_{i=1}^d u_i g'(X_i) \prod_{j\neq i} g(X_j)}{\prod_j g(X_j)} = \sum_{i=1}^d u_i \frac{g'(X_i)}{g(X_i)}.
\end{equation*}
Note that the above is a sum of independent random variables.
Moreover, these random variables have mean zero due to symmetry of $g$, and so by standard concentration inequalities, the value of the above sum is extremely likely to be on the order of its standard deviation.
Due to independence, the variance of the sum is the sum of the variances, and so
\begin{equation*}
	\var[\bigg]{\frac{\partial_u P(X)}{P(X)}} = \sum_{i=1}^d u_i^2 \var[\bigg]{\frac{g'(X_i)}{g(X_i)}} = a\norm{u}_2^2
\end{equation*}
for $a = \var{g'(Y)/g(Y)}$ where $Y$ is sampled from $(-C,C)$ according to $g$; this is, in particular, a constant (though not necessarily finite!).
We therefore see that, for sufficiently small $\delta$, the shift distance $\Delta(P,\delta u)$ should not depend on the dimension $d$.
Taking $\delta$ to be a small constant is equivalent to taking $C$ to be a large constant. 
This implies that such a density should exist so long as $\var{g'(Y)/g(Y)}$ is finite (which is in particular not the case for some candidate functions, such as the natural piecewise linear choice for $g$).

Despite the fact that the above outline can be formalised, our formal proof will take a different angle, though it can be understood through a very similar lens to the above.

\begin{lemma}
\label{lem:base}
	Fix a constant $C > 0$ and let $K = (-C,C)^d$ be a hypercube.
	There is a function $P\from \RR^d \to \RR$ supported on $K$ with $\int_{\RR^d} P(x) \d x = 1$ (i.e.\ $P$ is a probability density), $\mu(P) = 0$, and, for all $u\in \RR^d$ a unit vector, $\Delta(P,u) \leq \sqrt{3}/C$.
\end{lemma}

\begin{proof}
	We can define the desired function $P$ as follows.
	\begin{equation*}
		g(t) \defined \sqrt{\frac{3}{2C}} \max\set[\Big]{1 - \frac{\abs{t}}{C}, 0}, \qquad f(x) \defined \prod_{i=1}^d g(x_i), \qquad\text{and}\qquad P(x) \defined f(x)^2.
	\end{equation*}
	The following three facts follow from simple computations:
	\begin{equation*}
		\int g^2(t) \d t = 1, \qquad \int (g')^2(t) \d t = \frac{3}{C^2}, \qquad\text{and}\qquad \int g(t) g'(t) \d t = 0.
	\end{equation*}
	The first of the above tells us that $\int_{\RR^d} P(x) \d x = 1$, and so $P$ is indeed a probability density. Due to symmetry of $g$, we have that $\mu(P) = 0$.

	We will control $\Delta(P, u)$ in terms of $\norm{\partial_u f}_2^2 = \int_{\RR^d} \p[\big]{\partial_u f(x)}^2 \d x$.
	We may note that
	\begin{equation*}
		\partial_u f = \sum_i u_i \frac{\partial f}{\partial x_i} = \sum_i u_i g'(x_i) \prod_{j\neq i} g(x_j).
	\end{equation*}
	In expanding $(\partial_u f(x))^2$, it turns out that the integral of any cross term $(\partial f / \partial x_i)(\partial f / \partial x_j)$ with $i\neq j$ will be zero:
	\begin{equation*}
		\int_{\RR^d} \frac{\partial f}{\partial x_i} \, \frac{\partial f}{\partial x_j} \d x = \int_{\RR^d} g(x_i) g'(x_i) g(x_j) g'(x_j) \prod_{\l\notin\set{i,j}} g(x_\l)^2 \d x,
	\end{equation*}
	and we can pull out a factor of $\int_\RR g(t) g'(t) \d t = 0$, so the above is indeed zero. 
	We may therefore expand out $\norm[\big]{\partial_u f}_2^2$, disregarding any cross terms.
	\begin{equation*}
		\norm[\big]{\partial_u f}_2^2 = \int_{\RR^d}\p[\bigg]{\sum_{i=1}^d u_i \, \frac{\partial f}{\partial x_i}}^2 \d x = \int_{\RR^d}\p[\bigg]{\sum_i u_i^2 \p[\Big]{\frac{\partial f}{\partial x_i}}^2}\d x = \frac{3}{C^2}\norm{u}_2^2 = \frac{3}{C^2}.
	\end{equation*}
	We now turn our attention to the shift distance, and by expanding the definition find
	\begin{align*}
		\Delta(P,u) 
		&= \frac{1}{2}\int\abs[\big]{f^2(x + u) - f^2(x)} \d x\\
		&= \frac{1}{2} \int \abs[\big]{f(x + u) + f(x)} \, \abs[\big]{f(x + u) - f(x)} \d x\\
		&\leq \frac{1}{2} \norm[\big]{f(\blank + u) + f(\blank)}_2 \norm[\big]{f(\blank + u) - f(\blank)}_2,
	\end{align*}
	where the final inequality follows from Cauchy--Schwarz.
	By the triangle inequality, we have $\norm{f(\blank + u) + f(\blank)}_2 \leq 2$, so it remains to control $\norm{f(\blank + u) - f(\blank)}_2$.
	To do this, we may first note that $\norm{f(\blank + u) - f(\blank)}_2 \leq \norm{\partial_u f}_2$, as follows.
	\begin{align*}
		\norm[\big]{f(\blank + u) - f(\blank)}_2^2
		&= \int_{\RR^d} \p[\bigg]{\int_0^1 \partial_u f (x + tu) \d t}^2 \d x \\
		&\leq \int_{\RR^d} \int_0^1 \p[\big]{\partial_u f (x + tu)}^2 \d t \d x \\
		&= \int_0^1 \int_{\RR^d} \p[\big]{\partial_u f (x + tu)}^2 \d x \d t \\
		&= \norm[\big]{\partial_u f}_2^2.
	\end{align*}
	The first line above follows from the fundamental theorem of calculus, the second from Cauchy--Schwarz, and the third from Tonelli's theorem.
	We therefore have that, for any unit vector $u\in\RR^d$,
	\begin{equation*}
		\Delta(P, u) \leq \norm[\big]{\partial_u f}_2 = \frac{\sqrt 3}{C}.
	\end{equation*}
	This completes the proof of \Cref{lem:base}.
\end{proof}


\subsubsection{Proof of the Koml\'os conjecture}
\label{subsubsec:proof}

We now have all the pieces we need to put together a proof of the Koml\'os conjecture.

\begin{proof}[Proof of \Cref{thm:main}]
	While \Cref{thm:main} is stated for vectors with $\norm{v_i}_2 \leq 1$, we first show that it suffices to consider the case wherein $\norm{v_i}_2 = 1$ for all $i$.
	We will move from $\RR^d$ to $\RR^{d+n}$; let $\iota\from\RR^d\to\RR^{d+n}$ be the embedding map which pads a vector with $n$ zeroes, and $e_1,\dotsc,e_n$ be a basis of these added dimensions.
	We may then replace $v_i$ with $\iota(v_i) + \sqrt{1 - \norm{v_i}_2^2} \, e_i$.

	As the first $d$ coordinates of the vectors are unchanged, this cannot decrease the discrepancy.
	Moreover, as only a single vector has a component in the direction of each new basis vector $e_i$, this will not increase the discrepancy provided that the discrepancy is at least 1.
	In this way, we may assume that the Euclidean norm of $v_i$ is exactly 1 for all $i$; we may have significantly increased the dimension of the problem, but this does not matter in the rest of the proof.

	Fix $a=3$ and $c = 3\sqrt{3}$, let $K_0 \defined (-2ac,2ac)^d$, and let $P_0$ be the density on $K_0$ defined as the dilation (by a factor of $2a$) of the density in \Cref{lem:base}.
	We thus know that $P_0$ has $\Delta(P_0,2av) \leq \sqrt{3}/c$ for all $v\in \SS^{d-1} \sseq \RR^d$.
	Note that $(K_0,P_0)$ is a nice pair.

	We are given $v_1,\dotsc,v_n \in \RR^{d}$ with all $\norm{v_i}_2 = 1$ and must show that there are signs $\eps_1,\dotsc,\eps_n \in \signs$ such that $\sum_i \eps_i v_i \in K_0$.
	We will construct a probability distribution supported on $K_0$ with mean $\sum_i \eps_i v_i$.

	First, we inductively define the nice pair $(K_k,P_k)$ to be the $av_k$-lift of $(K_{k-1}, P_{k-1})$ for each $k=1,2,\dotsc,n$ in turn.
	Note that $K_k \sseq \RR^d \times \cube^k$ for all $k$, and that, by \Cref{lem:shift-distance-decreasing}, we have that $\Delta(P_k,u) \leq \Delta(P_0, u)$ for all $k$ and all $u\in \RR^d$.
	By \Cref{lem:lift-mass-discrete-fibre,lem:lift-mass-upstairs}, we find that 
	\begin{equation*}
		P_n\b[\big]{\upstairs{n}{k}} = \frac{1 - \Delta(P_{k-1}, 2av_k)}{2} \geq \frac{1 - \sqrt{3}/c}{2}.
	\end{equation*}
	In other words, the total mass of $P_n$ in any upper facet is at least some constant, and this constant can be made arbitrarily close to $1/2$ by increasing the value of $c$.

	We now iteratively apply collapses to this density: we will construct probability densities $Q_n,Q_{n-1},\dotsc,Q_1,Q_0$ on $K_n,\dotsc,K_0$ respectively, satisfying
	\begin{equation}
		\label{eq:mean-moving}
		\mu(Q_{k-1}) = \mu(Q_k) + \eps_k v_k.
	\end{equation}
	In particular, we will set $Q_n = P_n$, and then $Q_{k-1}$ will be an $f_k$-collapse of $Q_k$ for a suitable $av_k$-collapse function $f_k\from K_k \to K_{k-1}$.
	It is then immediate that each $Q_k$ is a density on $K_k$, and \Cref{lem:collapse-preserves-mass} tells us that, for all $\l\leq k\leq n$, we have
	\begin{equation*}
		Q_k\b[\big]{\upstairs{k}{\l}} \geq \frac{1 - \sqrt{3}/c}{2}.
	\end{equation*}
	It therefore suffices to prove that we can select the $av_k$-collapse function $f_k$ so that equality \eqref{eq:mean-moving} holds.
	But this follows from \Cref{lem:step}; for condition \eqref{eq:collapse-mass-upstairs} to hold, we require that 
	\begin{equation*}
		\frac{1 - \sqrt{3}/c}{2} \geq \frac{1}{2}\max\set[\Big]{1 - \frac{1}{a}, \, \frac{2}{a}},
	\end{equation*}
	and then we can conclude that our desired $av_k$-collapse function exists.
	After some simple rearrangements, noting that, as $a\geq 3$, we have $1 - 1/a \geq 2/a$, we conclude that $c \geq \sqrt{3}\,a$ suffices.
	Recalling that we set $a=3$, we see that $K_0 = (-18\sqrt{3},18\sqrt{3})^d$ suffices for the above argument to hold.
	In particular, $\mu(Q_n) = \mu(P_n) = \mu(P_0) = 0$, and so
	\begin{equation*}
		\mu(Q_0) = \sum_{i=1}^n \eps_i v_i,
	\end{equation*}
	and so the Koml\'os conjecture holds.
\end{proof}


\section{Counterexample to the Strong Koml\'os Conjecture}
\label{sec:strong-counter}

Now that the Koml\'os conjecture is known to be true, the next natural question is how far the solution can be pushed. 
Beyond just improving the value of the constant $C$ in \Cref{conj:komlos}, one can also ask whether a stronger version of the conjecture could also allow for a constant bound: can we bound all prefixes as well as the final sum?
What is the supremal value of
\begin{equation}
	\label{eq:strong-bound}
	\min_{\eps\in\signs^n} \; \max_{1\leq k \leq n} \; \norm[\bigg]{\sum_{i=1}^k \eps_i v_i}_\infty
\end{equation}
over all choices of sub-unit vectors $v_1,\dotsc,v_n \in \RR^d$?
Such questions date all the way back to Spencer's 1977 work \cite{Spe77} (i.e.\ before the Koml\'os conjecture was even posed), though the upper bound proved there was exponential in the dimension.
B\'ar\'any and Grinberg \cite{BG81} improved this to a bound linear in the dimension.
In terms of bounds which do not depend on the sequence length, results are very sparse. 
Indeed, the only other such upper bound which the author is aware of follows from OpenAI's work on the Euclidean Steinitz bound \cite{Clank26}, which proves an upper bound of $O(\sqrt{d})$.

In this section we present a proof of \Cref{thm:false}, discovered by Kintali \cite{Kin26}, showing that there is no constant upper bound to \eqref{eq:strong-bound}.
In fact, this will prove a lower bound of $\Omega(\sqrt{\log d})$.

We first present an outline and motivation of the proof and then give the formal details.


\subsection{Proof outline and motivation}

We now motivate the proof of \Cref{thm:false}, working backwards through the steps in the order they might be discovered.

\begin{figure}[htb]
    \centering
\begin{tikzpicture}[
  x=1cm,y=1cm,font=\normalsize,
  line cap=butt,line join=miter,
  axis/.style={draw=black,line width=.65pt,-{Stealth[length=5pt]}},
  boundary/.style={draw=black,line width=.8pt},
  tick/.style={draw=black,line width=.65pt},
  squareguide/.style={draw=black!25,densely dotted,line width=.45pt,
    line cap=round},
  region/.style={pattern=north east lines,pattern color=black!60},
  ticklabel/.style={font=\small,inner sep=1pt,fill=white,
    text height=2ex,text depth=.3ex},
  leftticklabel/.style={ticklabel,text height=1ex,text depth=0pt}]
  \def\cornerA{.75}
  \def\viewExtent{2.20}
  \def\rayExtent{4.00}
  \pgfmathsetmacro{\rotatedA}{sqrt(2)*\cornerA}

  \draw[squareguide] (-\cornerA,-\cornerA)
    rectangle (\cornerA,\cornerA);
  \begin{scope}
    \clip (-\viewExtent,-\viewExtent)
      rectangle (\viewExtent,\viewExtent);
    \foreach \sx in {-1,1} {
      \foreach \sy in {-1,1} {
        \fill[region] ({\sx*\cornerA},{\sy*\cornerA})
          rectangle ({\sx*\rayExtent},{\sy*\rayExtent});
        \draw[boundary] ({\sx*\rayExtent},{\sy*\cornerA})
          -- ({\sx*\cornerA},{\sy*\cornerA})
          -- ({\sx*\cornerA},{\sy*\rayExtent});
      }
    }
  \end{scope}
  \draw[axis] (-2.35,0) -- (2.35,0);
  \draw[axis] (0,-2.35) -- (0,2.35);
  \foreach \sgn in {-1,1} {
    \draw[tick] ({\sgn*\cornerA},-.055)
      -- ({\sgn*\cornerA},.055);
    \draw[tick] (-.055,{\sgn*\cornerA})
      -- (.055,{\sgn*\cornerA});
  }
  \node[leftticklabel] at (-\cornerA,.30) {$-a$};
  \node[leftticklabel] at (\cornerA,-.30) {$a$};
  \node[leftticklabel] at (.30,\cornerA) {$a$};
  \node[leftticklabel] at (-.30,-\cornerA) {$-a$};

  \draw[axis] (2.85,0) -- (5.05,0)
    node[midway,above=4pt,align=center] {rotate by\\[2pt]$\pi/4$};

  \begin{scope}[xshift=7.90cm]
    \draw[squareguide] (-\rotatedA,-\rotatedA)
      rectangle (\rotatedA,\rotatedA);
    \begin{scope}
      \clip (-\viewExtent,-\viewExtent)
        rectangle (\viewExtent,\viewExtent);
      \begin{scope}[rotate=45]
        \foreach \sx in {-1,1} {
          \foreach \sy in {-1,1} {
            \fill[region] ({\sx*\cornerA},{\sy*\cornerA})
              rectangle ({\sx*\rayExtent},{\sy*\rayExtent});
            \draw[boundary] ({\sx*\rayExtent},{\sy*\cornerA})
              -- ({\sx*\cornerA},{\sy*\cornerA})
              -- ({\sx*\cornerA},{\sy*\rayExtent});
          }
        }
      \end{scope}
    \end{scope}
    \draw[axis] (-2.35,0) -- (2.35,0);
    \draw[axis] (0,-2.35) -- (0,2.35);
    \foreach \sgn in {-1,1} {
      \draw[tick] ({\sgn*\rotatedA},-.055)
        -- ({\sgn*\rotatedA},.055);
      \draw[tick] (-.055,{\sgn*\rotatedA})
        -- (.055,{\sgn*\rotatedA});
    }
    \node[ticklabel] at ({-\rotatedA+.15},.55)
      {$-\sqrt{2}\,a$};
    \node[ticklabel] at ({\rotatedA-.15},-.55)
      {$\sqrt{2}\,a$};
    \node[ticklabel] at (.55,{\rotatedA-.15})
      {$\sqrt{2}\,a$};
    \node[ticklabel] at (-.55,{-\rotatedA+.15})
      {$-\sqrt{2}\,a$};
  \end{scope}
\end{tikzpicture}
	\caption{A diagram of how we can rotate coordinates in a two-dimensional space to go from two large entries to one larger entry.}
	\label{fig:rotate}
\end{figure}
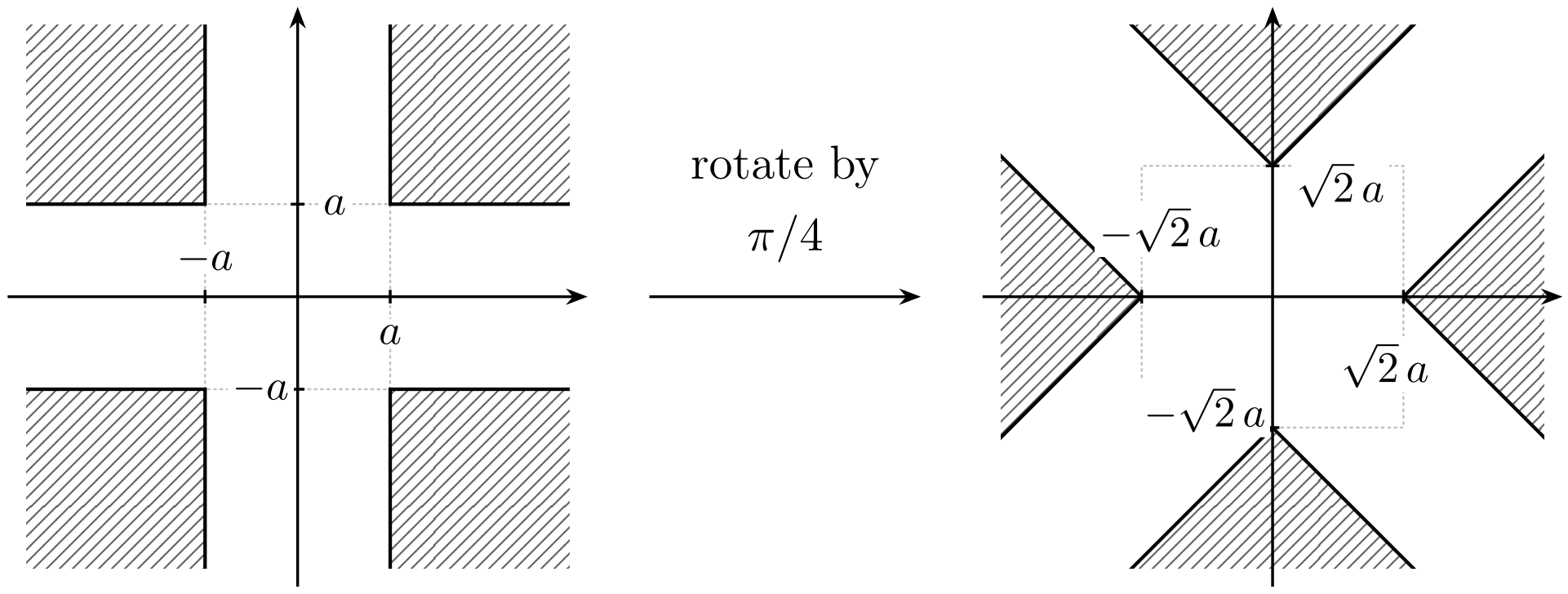

The first question is how to step from a sequence having a partial sum maximum entry at least $a$ in absolute value to an increased bound.
The idea behind this is that, if we can somehow force the sequence to simultaneously have two particular entries large (with our choice of which entries), then we can apply a rotation by $\pi/4$ to make a single entry larger; see \Cref{fig:rotate} for a diagram of this.

For ease of discussion of these results, we introduce some terminology now: say that a sequence $v_1,\dotsc,v_n$ of vectors in $\RR^d$ \emph{forces height $a$} if, for all vectors $p\in\RR^d$ (with no condition on the size of $p$) and all sequences of signs $\eps_1,\dotsc,\eps_n$, there is some $0\leq k \leq n$ such that $\norm[\big]{p + \sum_{i=1}^k \eps_i v_i}_\infty \geq a$.

So how do we construct a sequence in which the scenario in \Cref{fig:rotate} occurs? 
The key here is to be able to force simultaneity: it would be convenient to split the set of coordinates in half, force different things to happen in each half, and then assert that these two events must coincide.
While this might at first not seem possible, as the two halves share a sign sequence $\eps_1,\dotsc,\eps_n$, and we cannot control where in the sequence the large values are obtained, there is a way around this obstruction.

First, to be able to put vector sequences one after another in series and preserve their desirable properties (namely, having some partial sum with some entry/entries large), instead of just insisting that we should be able to force a large entry, we insist that this should happen even upon addition of any starting vector $p$.
Then, given two sequences $U$ and $V$, one can produce a larger sequence by putting them one after another; as $V$ can handle the addition of any starting vector $p$, it coming after the sequence $U$ will not stop it from working.

In fact, we can also combine $U$ and $V$ in a more complex way: if these are sequences on $\RR^d$, we can move to $\RR^{2d}$ (noting that there are no constraints on the dimension of the space that we're working in), put $U$ in the first $d$ coordinates of this larger space, and put $V$ in the last $d$ coordinates.
Then we can combine the sequences by inserting a copy of $V$ between every pair of entries of $U$, as well as at the start and end of the sequence; see \Cref{fig:too-many-columns} for an illustration of this.
Find where in $U$ the desired conclusion holds, and look at the copy of $V$ just after this point; as these two sequences operate on disjoint coordinates, they do not interfere with each other.
Therefore, $V$ can force some entry to be large while $U$ has already done this, allowing us to conclude that their implications happen simultaneously.

What do we want our sequences $U$ and $V$ to do? 
Each is on some number $D$ of coordinates, and to apply our rotation we would need to be able to conclude something of the form ``entries $i$ and $D+i$ are large simultaneously.''
To do this, we will have $U$ and $V$ perform different roles: $U$ will force many entries to be large simultaneously, and $V$ will force one entry large, with us having some control over which entry is large.

In particular, by carefully combining many copies of a sequence which forces one entry large, we will be able to produce a sequence which forces $d$ entries to be large simultaneously: this will be $U$.
Similarly, by inserting a contiguous copy of a sequence forcing one entry to be large in each set of $d$ coordinates of a larger space, we can force one entry of any given set of $d$ to be large.
Combining these two gives us exactly the conditions that we require.


\subsection{Formal proof}

We first introduce the key piece of terminology that we need in this proof.

\begin{definition}
	\label{def:forces-height}
	A sequence of vectors $v_1,\dotsc,v_n \in \RR^d$ is said to \emph{force height $a$} if for all $p\in \RR^d$ and all sign sequences $\eps_1,\dotsc,\eps_n \in\signs$ there is $0\leq k \leq n$ such that
	\begin{equation*}
		\norm[\bigg]{p + \sum_{i=1}^k \eps_i v_i}_\infty \geq a.
	\end{equation*}
	Moreover, say that the sequence of vectors \emph{$\l\times$forces height $a$} if for all $p$ and $\eps$ as above, there is some $k$ such that the vector $p + \sum_{i=1}^k \eps_i v_i$ has at least $\l$ different entries which are at least $a$ in absolute value.
\end{definition}

Note that we include the $k=0$ case in the above definition so that if $\norm{p}_\infty \geq a$ then we can take $p$ itself as the witness that the height of the sequence is at least $a$.

The key result, from which we will deduce \Cref{thm:false}, is the following.

\begin{proposition}
	\label{prop:strong-induction}
	If there is a sequence of $n$ sub-unit vectors in $\RR^d$ which forces height $a$ for some $a > 0$, then there is a sequence of at most $((n+1)d)^{2d}$ sub-unit vectors in $\RR^{2d^2}$ which forces height $\sqrt{2} a$.
\end{proposition}

Then, noting that the sequence $1$ of the unit vector in $\RR^1$ forces height $1/2$, it is immediate that there cannot be a constant bound on the maximum height of a sequence.
We perform the computations to check the bounds claimed in \Cref{thm:false} at the end of the section.

The proof of \Cref{prop:strong-induction} which we present has three steps, captured by \Cref{lem:d-large-values,lem:one-of-d,lem:lets-rotate-the-board}.

\begin{lemma}
	\label{lem:d-large-values}
	If the sequence $v_1,\dotsc,v_n \in \RR^d$ of sub-unit vectors forces height $a$, then there is an integer $N = (n+1)^d - 1$ and a sequence $u_1,\dotsc,u_N \in \RR^{d^2}$ of sub-unit vectors that $d\times$forces height $a$.
\end{lemma}

\Cref{lem:d-large-values} gives us half of the desired construction: we can force that at least $d$ entries of a partial sum are large simultaneously.

\begin{lemma}
	\label{lem:one-of-d}
	If the sequence $v_1,\dotsc,v_n \in \RR^d$ of sub-unit vectors forces height $a$ and $D \geq d$, then there is an integer $M=\binom{D}{d}n$ and a sequence $w_1,\dotsc,w_M \in \RR^D$ of sub-unit vectors such that, for all $p\in \RR^D$, all sequences $\eps_1,\dotsc,\eps_M \in \signs$ of signs, and all sets $X \sseq [D]$ of size $d$, there is $0 \leq k \leq M$ such that
	\begin{equation}
		\label{eq:force-one-of-d}
		\exists\,\l\in X \text{ such that } \abs[\bigg]{\p[\bigg]{p + \sum_{i=1}^k \eps_i w_i}_\l} \geq a.
	\end{equation}
\end{lemma}

\Cref{lem:one-of-d} is what gives us the second half of the target construction: we can force a large entry, but with some control over where: given any $d$ coordinates, we can force that one of those entries is large.

\begin{lemma}
	\label{lem:lets-rotate-the-board}
	If the sequence $v_1,\dotsc,v_n\in \RR^{2D}$ of sub-unit vectors is such that, for all $p\in \RR^{2D}$ and all sequences $\eps_1,\dotsc,\eps_n\in\signs$ of signs, there is $0 \leq k \leq n$ and $\l \in [D]$ such that
	\begin{equation}
		\label{eq:paired-forcing}
		\abs[\bigg]{\p[\bigg]{p + \sum_{i=1}^k \eps_i v_i}_\l} \geq a \quad \text{and}\quad \abs[\bigg]{\p[\bigg]{p + \sum_{i=1}^k \eps_i v_i}_{D + \l}} \geq a,
	\end{equation}
	then there is a sequence $u_1,\dotsc,u_n \in \RR^{2D}$ of sub-unit vectors which forces height $\sqrt{2}a$.
\end{lemma}

\Cref{lem:lets-rotate-the-board} actually performs the desired boost: if we can produce a sequence of vectors in $\RR^{2D}$ and a partition of the coordinates into pairs such that there is some partial sum in which both entries in some pair are large simultaneously, then we can convert this into a sequence which forces height $\sqrt{2}a$ through the ``rotation'' step discussed in the proof motivation.

We first show how to combine the above three lemmas into a proof of \Cref{prop:strong-induction} and then prove the three lemmas in turn.
The idea behind the proof is that we can insert a copy of the sequence produced by \Cref{lem:one-of-d} between every pair of vectors in the sequence produced by \Cref{lem:d-large-values}, where these two sequences are given disjoint coordinate sets, as shown in \Cref{fig:too-many-columns}.
This is enough to force the pairing property required by \Cref{lem:lets-rotate-the-board} and so give us the required boost.

\begin{figure}[htb]
    \centering
	\begin{equation*}\boxed{\begin{matrix}
		  &   &  &  & \uparrow &  &  &  & \uparrow &  &  &  & \uparrow &   &  &   & \uparrow &  &  &  \\
		0 & 0 & \dotsc & 0 & u_1 & 0 & \dotsc & 0 & u_2 & 0 & \dotsc & 0 & u_3 &   & \boldsymbol{\cdot}\boldsymbol{\cdot}\boldsymbol{\cdot} &   & u_N & 0 & \dotsc & 0 \\
		  &   &  &  & \downarrow &  &  &  & \downarrow &  &  &  & \downarrow &   &  &   & \downarrow &  &  &  \\
		\uparrow & \uparrow &  & \uparrow &  & \uparrow &  & \uparrow &  & \uparrow &  & \uparrow &  &   &  &   &  & \uparrow &  & \uparrow \\
		w_1 & w_2 & \dotsc & w_M & 0 & w_1 & \dotsc & w_M & 0 & w_1 & \dotsc & w_M & 0 &  & \boldsymbol{\cdot}\boldsymbol{\cdot}\boldsymbol{\cdot} &  & 0 & w_1 & \dotsc & w_M \\
		\downarrow & \downarrow &  & \downarrow &  & \downarrow &  & \downarrow &  & \downarrow &  & \downarrow &  &   &  &   &  & \downarrow &  & \downarrow \\
	\end{matrix}}\end{equation*}
	\caption{The constructed matrix in the proof of \Cref{prop:strong-induction}}
	\label{fig:too-many-columns}
\end{figure}

\begin{proof}[Proof of \Cref{prop:strong-induction}]
	Assume that $v_1,\dotsc,v_n \in \RR^d$ is a sequence of sub-unit vectors which forces height $a$, for some $a > 0$.
	Apply \Cref{lem:d-large-values} to produce a sequence $u_1,\dotsc,u_N \in \RR^{d^2}$ which forces $d$ entries at least $a$ simultaneously, and apply \Cref{lem:one-of-d} with $D = d^2$ to produce a sequence $w_1,\dotsc,w_M \in \RR^{d^2}$ which forces one of any particular $d$ entries to be at least $a$.
	Note that $M = \binom{D}{d}n \leq d^{2d} n$.

	We combine these two sequences into a single sequence of vectors in $\RR^{2d^2}$ as follows.
	First embed the sequences into $\RR^{2d^2}$ by putting each vector $u_1,\dotsc,u_N$ in coordinates $1,\dotsc,d^2$ and filling the rest with zeroes, and conversely putting each vector $w_1,\dotsc,w_M$ in coordinates $d^2+1,\dotsc,2d^2$ and filling the first half with zeroes.
	We slightly abuse notation by continuing to refer to these vectors as $u_1,\dotsc,u_N$ and $w_1,\dotsc,w_M$.
	
	Insert a copy of the sequence $w_1,\dotsc,w_M$ between every $u_i$ and $u_{i+1}$, as well as before $u_1$ and after $u_N$, as shown in \Cref{fig:too-many-columns}.
	Call the resulting sequence $x_1,\dotsc,x_L$, where 
	\begin{equation*}
		L = N + (N+1)M \leq (n+1)^d\p[\big]{nd^{2d}+1} \leq \p[\big]{(n+1)d}^{2d},
	\end{equation*}
	as claimed.

	Note that every vector $x_i$ has all its non-zero entries in either the first half or the second half of the coordinates; we can thus consider the entries on either half separately.
	
	Firstly, the (now very spread out) copy of the sequence $u_1,\dotsc,u_N$ forces $d$ different entries to be at least $a$ in absolute value simultaneously. 
	Let $k$ be the index in the sequence $x$ where this occurs---that is, for all $p\in\RR^{2d^2}$ and $\eps\in\signs^L$, there is a set $X\sseq[d^2]$ with $\abs{X} = d$ and $0\leq k\leq L$ such that 
	\begin{equation*}
		\forall\;\l \in X, \; \abs[\bigg]{\,\p[\bigg]{p + \sum_{i=1}^k \eps_i x_i}_\l} \geq a.
	\end{equation*}
	Let $Y$ be $X$ translated to the second half of the coordinates, i.e.\ $Y \defined \set{\l + d^2 \st \l \in X}$.
	Consider the copy of $w_1,\dotsc,w_M$ occurring just after the element of the sequence $u_i$ corresponding to position $k$ (or the copy before $u_1$ if $k < M$).
	Property \eqref{eq:force-one-of-d} of $w_1,\dotsc,w_M$ implies that this sequence forces some entry in $Y$ to be at least $a$.

	This implies that the sequence $x_1,\dotsc,x_L$ forces entries $\l$ and $d^2 + \l$ to simultaneously be at least $a$ in absolute value at some point, which is precisely as required by condition \eqref{eq:paired-forcing} of \Cref{lem:lets-rotate-the-board}.
	This means that we may apply \Cref{lem:lets-rotate-the-board} to find a sequence which forces height $\sqrt{2} a$, as required, proving \Cref{prop:strong-induction}.
\end{proof}

It remains to prove the three lemmas.

\medskip

The proof of \Cref{lem:d-large-values} is inductive: to go from $\l\times$forcing to $(\l+1)\times$forcing, we can insert a copy of our original forcing sequence between every pair of elements, using $d$ new coordinates; this provides an extra large entry.

\begin{proof}[Proof of \Cref{lem:d-large-values}]
	We inductively prove that, for all $\l\geq 1$, there is a sequence $V_\l$ of $(n+1)^\l - 1$ vectors in $\RR^{d\l}$ which $\l\times$forces height $a$.

	For the base case, let $V_1 = (v_1,\dotsc,v_n)$ be the sequence given in the lemma. 
	Then assume for induction that we have constructed a sequence $V_\l = (u_1,\dotsc,u_m)$ of vectors in $\RR^{d\l}$ for some $\l\geq 1$ which $\l\times$forces height $a$.
	We construct $V_{\l + 1}$ by taking $V_\l$ and adding $d$ new entries to each vector, all zero.
	Then insert a copy of the original sequence $v_1,\dotsc,v_n$ living in the $d$ new coordinates between $u_i$ and $u_{i+1}$ for each $i$, as well as before $u_1$ and after $u_m$.
	This is essentially the same as the construction shown in \Cref{fig:too-many-columns}.

	We know by inductive hypothesis that there is some position $k$ such that the sequence $V_\l$ has at least $\l$ different entries large simultaneously.
	Look at the copy of $v_1,\dotsc,v_n$ just after this point: it itself forces some entry large, and lives in a disjoint set of coordinates to $V_\l$.
	Therefore, when it forces some coordinate to be large, the $\l$ large coordinates from $V_\l$ are also still large, so we get $\l + 1$ large coordinates in total.
	Note moreover that the length condition also carries through: the number of vectors in $V_{\l + 1}$ is
	\begin{equation*}
		m + (m+1)n = \p[\big]{(n+1)^\l - 1} + (n+1)^\l n = (n+1)^{\l + 1} - 1,
	\end{equation*}
	as required.
\end{proof}

To prove \Cref{lem:one-of-d}, we can simply put a contiguous copy of $v_1,\dotsc,v_n$ into each possible set of $d$ coordinates in the larger space $\RR^D$.

\begin{proof}[Proof of \Cref{lem:one-of-d}]
	We directly construct the required sequence out of copies of $v_1,\dotsc,v_n$.
	Enumerate the sets $X\sseq [D]$ of size $d$ as $X_1,\dotsc,X_m$ for $m = \binom{D}{d}$.
	Given $i\leq m$, let $V_i$ be the copy of the sequence $v_1,\dotsc,v_n$ placed in $\RR^D$ using the coordinates indexed by $X_i$ and setting other entries to zero.
	Then let $w_1,\dotsc,w_M$ be formed by appending the sequences $V_1,V_2,\dotsc,V_m$ one after another.

	Given any particular $X$, there is a contiguous copy of $v_1,\dotsc,v_n$ in the sequence $w_1,\dotsc,w_M$ living on the coordinates in $X$.
	This forces height $a$, and the entry forced to be large is in $X$, as required.
	Note that $M = \binom{D}{d} n$, as claimed.
\end{proof}

For the final lemma, \Cref{lem:lets-rotate-the-board}, we simply apply a suitable rotation of space, as shown in \Cref{fig:rotate}.

\begin{proof}[Proof of \Cref{lem:lets-rotate-the-board}]
	Let $P$ be the following $2\times 2$ real matrix.
	\begin{equation*}
		P = \p[\bigg]{\begin{matrix} 1/\sqrt{2} & -1/\sqrt{2} \\ 1/\sqrt{2} & 1/\sqrt{2} \end{matrix}}.
	\end{equation*} 
	Note that $P$ applies a rotation of $\pi / 4$ to the plane.
	Let $e_1,\dotsc,e_{2D}$ be the standard basis for $\RR^{2D}$.
	Then let $Q$ be the $2D \times 2D$ block matrix consisting of $D$ copies of $P$, one copy to rotate the plane spanned by $e_\l$ and $e_{D+\l}$ for all $1\leq\l \leq D$.
	We define $u_i = Qv_i$, and claim this suffices.
	
	Indeed, assume that condition \eqref{eq:paired-forcing} holds for some $\l$ for starting vector $Q^{-1} p$.
	Note that, if vector $x\in\RR^2$ has $\abs{x_1}\geq a$ and $\abs{x_2}\geq a$, then
	\begin{equation*}
		\norm[\big]{P x}_\infty \geq \sqrt{2} a.
	\end{equation*}
	This proves \Cref{lem:lets-rotate-the-board}.
\end{proof}

With the three lemmas proved, so too is \Cref{prop:strong-induction}.
It remains to check that the dependencies between the variables are as claimed.

\begin{proof}[Proof of \Cref{thm:false}]
	We iteratively apply \Cref{prop:strong-induction}, starting with the base of a single vector in $\RR^1$; this forces height $1/2$, and has $n=d=1$.
	Each application of \Cref{prop:strong-induction} sends $a$ to $\sqrt{2}a$, $d$ to $2d^2$, and $n$ to at most $((n+1)d)^{2d}$.
	Let $a_k$, $d_k$ and $n_k$ be the values of these parameters after $k$ of the above steps.
	It is clear that $a_k = 2^{k/2 - 1}$ and $d_k = 2^{\p{2^k-1}}$.
	We claim that
	\begin{equation*}
		n_k \leq 2^{2^{2^{k+1}}} = 2^{4d_k^2}.
	\end{equation*}
	Indeed, we may compute
	\begin{align*}
		\log_2 n_{k+1} 
		&\leq 2d_k \p[\big]{\log_2(n_k + 1) + \log_2 d_k} \\
		&\leq 2d_k \p[\big]{4d_k^2 +1 + \log_2 d_k} \\
		&\leq 16d_k^4 \\
		&= 4d_{k+1}^2,
	\end{align*}
	as claimed.

	Now, given an integer $d\geq 1$, take $k$ maximal such that $d_k \leq d$.
	Maximality of $k$ implies that $d_k \geq \sqrt{d/2}$.
	Note that
	\begin{equation*}
		n_k \leq 2^{4d_k^2} \leq 2^{4d^2} = \exp\p[\big]{O(d^2)},
	\end{equation*}
	as claimed in the statement of the theorem. 
	It remains to show that $a_k \geq c\sqrt{\log d}$.
	Recalling that $\log_2 d_k = 2^k - 1$ and that $a_k = 2^{k/2 - 1}$, we find that
	\begin{equation*}
		a_k = \frac{1}{2}\sqrt{1 + \log_2 d_k} \geq \sqrt{\tfrac{1}{4}\p[\big]{\tfrac{1}{2}\log_2 d + \tfrac{1}{2}}} \geq c\sqrt{\log d}
	\end{equation*}
	for some constant $c > 0$.
	Thus \Cref{thm:false} is proved.
\end{proof}


\section{Conclusion}
\label{sec:conc}

We have presented motivated proofs of the Koml\'os conjecture along with a counterexample to the strong Koml\'os conjecture.
It seems likely that the ideas contained in the proof will have further applications.
Indeed, some applications have already been found, as the technique from the original proof \cite{GFL26a} was used by OpenAI's internal model to determine the Euclidean Steinitz bound \cite{Clank26}.
It would be interesting to know whether OpenAI's (quite technical) proof in \cite{Clank26} could likewise be made elementary.

This does raise the question of what problems remain open.
Inspired by Banaszczyk's 2012 proof \cite{Ban12} that the Steinitz bound on the maximal prefix $\norm[\big]{\sum_{i=1}^k \eps_i v_i}_2$ is at most $O(\sqrt{d})$ if one is also allowed to permute the order of the vectors, we pose the following question.

\begin{question}
	Is there a constant $C$ such that the following holds?
	Given a sequence $v_1,\dotsc,v_n\in\RR^d$ of vectors satisfying $\norm{v_i}_2 \leq 1$ for all $i\leq n$, there is a permutation $\sigma\in S_n$ and sequence of signs $\eps_1,\dotsc,\eps_n\in\signs$ such that
	\begin{equation*}
		\norm[\bigg]{\sum_{i=1}^k \eps_i v_{\sigma(i)}}_\infty \leq C \quad \text{ for all } \; k \leq n.
	\end{equation*}
\end{question}

We remark that, as noted by Kintali \cite{Kin26}, the construction described in \Cref{sec:strong-counter} does not necessarily provide a counterexample to the above as many vectors are repeated; we can rearrange to have all the repeated vectors appear contiguously and give them signs so that they cancel out.

Another natural question is the probabilistic version of the Koml\'os conjecture: if $\eps_i\in\signs$ are chosen independently and uniformly at random, how small can $\prob{\norm{\sum_i \eps_i v_i}_\infty \leq C}$ be (as a function of $C$)?
A lower bound on this probability is $\exp(-O(n/C^2))$ for sufficiently large constant $C$, due to recent work of Akbas and Sra \cite{AS26a}. 
A matching upper bound of $\exp(-\Omega(n/C^2))$ can be found by taking copies of basis vectors, provided that $C$ and $n$ are sufficiently large ($n$ large in terms of $C$), though the value of the constant in the exponent remains unknown.

\begin{question}
	Let $\xi_1,\dotsc,\xi_n$ be independent Rademacher random variables, each uniformly distributed on the set $\signs$.
	For positive integer $n$ and constant $C > 0$, define
	\begin{equation*}
		p_n(C) \defined \inf_{d,v_1,\dotsc,v_n} \prob[\bigg]{\norm[\Big]{\sum_{i=1}^n \xi_i v_i}_\infty \leq C},
	\end{equation*}
	where the infimum is over all integers $d\geq 1$ and sub-unit vectors $v_1,\dotsc,v_n\in\RR^d$.
	What is the value of the limit
	\begin{equation*}
		\lim_{C\to\infty} \; \lim_{n\to\infty} \frac{-C^2}{n} \log p_n(C),
	\end{equation*}
	provided that it exists?
\end{question}

Despite the strong Koml\'os conjecture being false, it remains open what the true growth rate of the maximal infinity norm of a prefix sum is.

\begin{question}
	Fix $a > 0$.
	What is the minimum $d > 0$ such that there exists $n$ and a sequence $v_1,\dotsc,v_n\in\RR^d$ of sub-unit vectors such that, for every choice $\eps_1,\dotsc,\eps_n\in\signs$, there exists $k\leq n$ with $\norm[\big]{\sum_{i=1}^k \eps_i v_i}_\infty \geq a$?
	More precisely, define
	\begin{equation*}
		H(d) \defined \sup_n \sup_{\substack{v_1,\dotsc,v_n\in\RR^d \\ \norm{v_i}_2 \leq 1}} \; \min_{\eps\in\signs^n} \max_{k\leq n} \; \norm[\bigg]{\sum_{i=1}^k \eps_i v_i}_\infty.
	\end{equation*}
	What is the value of $H(d)$?
\end{question}
The present best-known bounds on $H(d)$ are $\Omega(\sqrt{\log d})$, as presented in \Cref{thm:false} and due to Kintali \cite{Kin26}, and $O(\sqrt{d})$, due to the recent OpenAI proof \cite{Clank26} of the Euclidean Steinitz result.

Finally, and perhaps most ambitiously, one could ask for the precise minimum value of $C$ for which the Koml\'os conjecture holds; the present best lower bound is $1+\sqrt{2}$ due to Kunisky \cite{Kun21}, and the best claimed upper bound appears to be approximately $6.9$ due to Akbas and Sra \cite{AS26b}.
This problem seems particularly out of reach and is, in the author's opinion, likely to stand unsolved for at least several more weeks.


\bibliographystyle{abbrvnat}  
\renewcommand{\bibname}{Bibliography}
\bibliography{main}


\end{document}